\documentclass{article}

\usepackage{tikz}
\usepackage{amsmath}
\usepackage{amsfonts}
\usepackage{amssymb}
\usepackage{xspace}
\usepackage{graphicx}
\usepackage[title]{appendix}
\usepackage{color}
\usepackage{epsf,epsfig}
\usepackage{graphics,color}
\usepackage{amsmath}
\usepackage{amssymb}
\usepackage{cite}
\usepackage{verbatim}
\usepackage{float}
\usepackage{graphicx}
\usepackage{amsthm}
\usepackage{textcomp}
\usepackage{subfig}
\usepackage{hyperref}
\usepackage{yhmath}
\usepackage{mathrsfs}
\numberwithin{equation}{section}
\makeatletter
\newcommand{\res}      {\mathop{\hbox{\vrule height 7pt width .5pt depth
			0pt\vrule height .5pt width 6pt depth 0pt}}\nolimits}

\newenvironment{prooof}{\par\noindent{\sc Proof:}
}{\hfill\llap{$\Box$}\vspace{1\baselineskip}\par\noindent}

\newenvironment{proofof}{\par\noindent{\sc Proof}
}{\hfill\llap{$\Box$}\vspace{1\baselineskip}\par\noindent}

\newcommand{\adl}{\vspace{1\baselineskip}}

\newtheorem{theorem}{Theorem}[section]
\newtheorem{proposition}[theorem]{Proposition}
\newtheorem{lemma}[theorem]{Lemma}
\newtheorem{corollary}[theorem]{Corollary}
\newtheorem{remark}[theorem]{Remark}
\newtheorem{definition}[theorem]{Definition}
\newtheorem{example}[theorem]{Example}

\newcommand{\beq}{\begin{equation}}
\newcommand{\eeq}{\end{equation}}
\newcommand{\rif}[1]{(\ref{#1})}
\newcommand{\ba}{\begin{array}}
\newcommand{\ea}{\end{array}}
\newcommand{\bt}{\begin{theorem}}
\newcommand{\et}{\end{theorem}}
\newcommand{\bp}{\begin{proposition}}
\newcommand{\ep}{\end{proposition}}
\newcommand{\bl}{\begin{lemma}}
\newcommand{\el}{\end{lemma}}
\newcommand{\bc}{\begin{corollary}}
\newcommand{\ec}{\end{corollary}}
\newcommand{\bi}{\begin{itemize}}
\newcommand{\ei}{\end{itemize}}
\newcommand{\ben}{\begin{enumerate}}
\newcommand{\een}{\end{enumerate}}
\newcommand{\bpf}{\begin{prooof}}
\newcommand{\epf}{\end{prooof}}
\newcommand{\bpff}{\begin{proofof}}
\newcommand{\epff}{\end{proofof}}
\newcommand{\bdf}{\begin{definition}\rm}
\newcommand{\edf}{\end{definition}}
\newcommand{\br}{\begin{remark}\rm}
\newcommand{\er}{\end{remark}}
\newcommand{\bex}{\begin{example}\rm}
\newcommand{\eex}{\end{example}}
\def\qu#1{{[\![\, #1\, ]\!]}}
\def\pri{\hbox to 10pt{\hfil\hbox to 0.4pt{\vrule height5pt width0.4pt
                 depth0pt}\vrule width5pt height0.4pt depth0pt\hfil}}
\newcommand{\jump}[1]{\text{{\rm \textlbrackdbl}}{#1}\text{{\rm \textrbrackdbl}}}

\newcommand{\gM}{{\bf M}}

\newcommand{\R}{{\cal R}}

\newcommand{\calL}{{\cal L}}

\newcommand{\gR}{{\mathbb R}}

\newcommand{\Nat}{{\mathbb N}}
\newcommand{\Rn}{{\mathbb R}^{n}}
\newcommand{\RN}{{\mathbb R}^{N}}

\newcommand{\RnN}{{\mathbb R}^{n+N}}

\newcommand{\Om}{\Omega}
\newcommand{\A}{{\cal A}}

\newcommand{\Ha}{{\cal H}}

\newcommand{\cart}{\mathop{\rm cart}\nolimits}
\def\det{\mathop{\rm det}\nolimits}

\newcommand{\If}{\ \mbox{\rm if}\ }

\newcommand{\wc}{\rightharpoonup}

\def\Lip{\mathop{\rm Lip}\nolimits}

\newcommand{\Sph}{{{\mathbb S}}}

\def\Det{\mathop{\rm Det}\nolimits}

\newcommand{\gn}{{\bf n}}

\newcommand{\gr}{{\bf r}}

\def\dd{\mathrm{d}}

\let\d=\delta
\let\e=\varepsilon
\let\f=\phi
\let\vf=\varphi
\let\g=\gamma

\let\l=\lambda
\let\m=\mu

\let\o=\omega
\let\p=\pi

\let\r=\rho
\let\s=\sigma
\let\t=\theta

\let\x=\xi
\let\y=\eta
\let\z=\zeta
\let\vf=\varphi

\let\LL=\Lambda
\let\SS=\Sigma
\let\O=\Omega

\let\vect=\overrightarrow
\let\wih=\widehat
\let\wdg=\wedge

\let\wid=\widetilde

\let\pa=\partial
\let\sb=\subset

\let\fa=\forall
\let\tim=\times
\let\mx=\mbox

\let\sm=\setminus
\let\ol=\overline
\let\ul=\underline
\let\ds=\displaystyle

\let\lan=\langle
\let\ran=\rangle
\let\i=\infty

\let\nolm=\nolimits

\let\pa=\partial
\let\sb=\subset

\let\fa=\forall
\let\tim=\times
\let\mx=\mbox

\let\sm=\setminus
\let\ol=\overline
\let\ul=\underline
\let\ds=\displaystyle

\let\ii=\infty

\def\LL{\mathop{\Lambda}}

\def\gn{{\mathbb{N}}}

\def\gr{{\mathbb{R}}}

\def\GM{{\mathbf{M}}}

\def\cA{{\mathcal{A}}}

\def\cD{{\mathcal{D}}}
\def\cE{{\mathcal{E}}}
\def\cF{{\mathcal{F}}}
\def\cG{{\mathcal{G}}}
\def\cH{{\mathcal{H}}}
\def\cI{{\mathcal{I}}}

\def\cK{{\mathcal{K}}}
\def\cL{{\mathcal{L}}}

\def\cR{{\mathcal{R}}}
\def\cS{{\mathcal{S}}}
\def\cT{{\mathcal{T}}}

\def\cV{{\mathcal{V}}}

\let\d=\delta
\let\e=\varepsilon
\let\f=\phi
\let\g=\gamma

\let\l=\lambda
\let\m=\mu

\let\o=\omega
\let\p=\pi

\let\r=\rho
\let\s=\sigma
\let\t=\theta

\let\x=\xi
\let\y=\eta
\let\z=\zeta

\let\vf=\varphi

\let\O=\Omega

\def\wih{\widehat}
\def\wid{\widetilde}

\newcommand{\rP}{{\mathbb P}}

\title{\Large \bf
Explicit formulas of strict BV relaxed energies \\ for polyconvex functionals with linear growth}
\author{\it Domenico Mucci
\footnote{
Dipartimento di Scienze Fisiche, Informatiche e Matematiche, Universit\`{a} degli studi di Modena e Reggio Emilia. Via Giuseppe Campi 213/b, I-41125 Modena (Italy). e-mail:
	domenico.mucci@unimore.it}
and Riccardo Scala
\footnote{Dipartimento di Ingegneria dell'Informazione e Scienze Matematiche, Universit\`a di Siena. Via Roma 56, I-53100 Siena (Italy).
e-mail: riccardo.scala@unisi.it}
}
\begin{document}
%\date{}
%\topskip=1.5truecm
\maketitle
\topskip=1.5truecm
       %%%%%%%%% Abstract %%%%%%%%%%%%%%%%
%
{\small{\bf Abstract.} We consider polyconvex functionals defined on vector valued functions, the model case being the graph area functional, and we analyze the relaxed energy with respect to the strict convergence in $BV$. In several cases where the relaxed area is known, by exploiting a continuity result by Reshetnyak we are able to find an explicit formula for wide classes of
integrands with linear growth in the minors of the gradient.
We preliminarily extend to the $BV$ setting a continuity property observed by Acerbi--Dal Maso in the Sobolev case.
Finally, partial results concerning the relaxation of the vortex map with respect to the $L^1$ convergence are obtained.}

\adl\par\noindent
\adl\par\noindent
%{46E35; 26B35; 49Q15%; 49Q20
%}
%
\bigskip
\par\noindent

\maketitle

\section{Introduction}

Since the seminal work of Ball \cite{Ball}, the concept of polyconvexity has played a fundamental role in continuum mechanics, particularly within the realm of nonlinear elasticity and related analytical-geometric problems. Under standard growth
{conditions} the involved energies turn out to be lower semicontinuous, whereas particular {attention} has been payed to the non standard cases, including the one when the functional has a linear growth. This is, for instance,  the case of the area and the total variation functionals. In \cite{AcDM}, the authors analyse continuity properties of polyconvex {energies} in the linear growth regime, proving that in general these are not lower semicontinuous and hence a relaxation procedure is needed, leading to lower semicontinuous functionals that however are not explicit. More precisely, if $f:\Om\times \gR^N\times \gR^{N\times n}{\to\gR}$ is a polyconvex integrand, with $\Om\subseteq\gR^n$ {a bounded open set,} $n,N\geq  2$, and $f$ has a linear growth of the form
\begin{equation}\label{eq1}
c\,|M(G)|\leq f(x,y,G)\leq C\,|M(G)|	\qquad \forall\, x\in \Om,\;y\in \gR^N,\;G\in \gR^{N\times n},
\end{equation}
with $c,C>0$,
one defines, for every $u\in L^1(\Om;\gR^N)$, the relaxed functional
 $$  \wid\cF_{L^1}(u,\Om)=\inf\left\{\liminf_{h\to\i}\int_{\Om} f(x,u_h(x),\nabla u_h(x))\dd x \mid \{u_h\}\sb C^1(\Om;\RN)\,,\,\,u_h\to u\,\,\text{in}\,\,L^1(\O;\RN)\right\} .
 $$
Above, $M(G)$ is the vector containing all the determinants of the $(k\times k)$-minors of $G$, $k=0,\dots,\min\{n,N\}$ (see Section \ref{sec:21} below). A special case of main interest is the area functional obtained by setting $f(x,y,G)=|M(G)|$, whose relaxation as above {is denoted} by $\wid\cA_{L^1}(u,\Om)$.
In \cite{AcDM} it has been shown {(under suitable assumptions on $f$)} that $\wid\cF_{L^1}(u,\Om)$ gives rise to a lower semicontinuous functional in $L^1(\Om;\gR^N)$ that extends the energy
 $$\mathcal E_f(u,\Om):=\int_\Om f(x,u(x),\nabla u(x))\dd x$$
for $C^1$ maps, and that however does not admit {an} integral representation; more precisely, the set function $A\mapsto {\wid\cA}_{L^1}(u,A)$, defined for all open set $A\subseteq \Om$, is not subadditive, and then it is not the restriction of any Borel measure. The latter assertion was originally conjectured by De Giorgi in \cite{DG}, who suggested to prove the lack of {subadditivity, in case $n=N=2$,} for the vortex map $u_V$ and the triple junction function $u_T$, which {become} the unique two nontrivial examples on which the relaxed area functional can be explicitly computed (see \cite{BES1} and \cite{Sc20}).

The issue regarding the lack of locality of the relaxed functional ${\wid\cA}_{L^1}(u,\Om)$ is closely linked to the wide freedom of approximation imposed by the $L^1$-norm. From the observation that the domain of $\wid\cF_{L^1}(u,\Om)$ is a subset of $BV(\Om;\gR^N)$,
%it is natural
{in order to overcome the previous drawbacks, a possible strategy is} to replace the $L^1$ convergence with the strict convergence in $BV$, leading to the new relaxed functional
\begin{equation*}
\wid\cF_{BV}(u,\Om)=\inf\left\{\liminf_{h\to\i}\mathcal E_f(u_h,\Om) \mid \{u_h\}\sb C^1(\Om;\RN)\,,\,\,u_h\to u\,\,\text{strictly in}\,\,BV(\Om;\RN)\right\}.
\end{equation*}
The corresponding relaxation of the area functional in the strict convergence in $BV$ is denoted by $\wid\cA_{BV}(u,\Om)$, {and the} analysis of the latter energy has been addressed in \cite{Mu22,BSS,BSS24,Ca24,CM}. In particular, integral representation formulas for  $\wid\cA_{BV}(u,\Om)$ has been {shown} in the case $n=2$ in \cite{BSS,BSS24,Ca24}, when for instance $u$ is a $0$-homogeneous or piecewise Lipschitz map, or even a Sobolev map with values in $\mathbb S^1$. The latter case has been generalized to the higher dimension $n\geq 3$ in \cite{CM}, where however it is proved that the set function $A\mapsto \wid\cA_{BV}(u,A)$ is still non-subadditive when $N\geq3$.

In the $2$-dimensional case $n=2$, a recent work \cite{Sc} addressed the proof of locality of $\wid\cF_{BV}(u,\Om)$, namely
{showing} that the set function $\Om\supseteq A\mapsto \wid\cF_{BV}(u,A)$ is actually the restriction to open sets of a Borel measure. In light of this result, hopes grow for finding a general integral representation for $\wid\cF_{BV}(u,\Om)$.
%
%The fact that in the case of strict convergence the relaxation \alert{procedure gives rise to a local functional (in dimension $n=2$), is due to the fact that} approximating a function $u\in BV(\Om;\gR^N)$ is more stringent than approximating it merely in $L^1$;

{The locality result obtained in \cite{Sc}
relies on the fact that in dimension $n=2$, and for any $N\geq 2$, the strict convergence in $BV$ is strong enough to ensure that suitable gluing techniques can be performed by paying a small amount of energy. Notice that this cannot happen in the case of $L^1$ convergence, where concentration effects may appear.}

{We expect that the locality property continues to hold in the codimension $N=2$ case, for any $n\geq 2$.}
More specifically, if $\{v_k\}$ is a sequence of smooth functions approximating $u$ in $L^1(\Om;\gR^N)$ and such that $ \mathcal E_f(v_k,\Om)<C<+\infty$ for all $k$, the first inequality in \eqref{eq1} implies that the graphs $G_{v_k}$ of $v_k$, seen as integer multiplicity {rectifiable currents in} $\Om\times \gR^N$, have uniformly bounded masses, and then, up to subsequences, they converge to some {Cartesian} current $T$ with underlying map $u$. {In general, there are a lot of such possible {limit currents} $T$ for graphs $G_{v_k}$ of qualitatively different $L^1$-approximating sequences, and this lack of uniqueness determines the non locality feature addressed in \cite{AcDM}.}
%in particular $T$ is not uniquely determined.

In contrast, {if $N=2$, for every sequence $v_k$ converging to $u$ strictly in $BV(\Om;\gR^2)$ and with equibounded area, the limit current $T$ is unique} (this is a result in \cite{Mu22}).  Hence, any function $u\in BV(\Om;\gR^2)$ {satisfying $\wid\cA_{BV}(u,\Om)<\infty$ has a unique corresponding Cartesian current $T_u$ (called {minimal completely vertical lifting} of $u$) which comes into the play in the relaxation approach.}
%%coincides with any limit point of $G_{v_k}$ for any sequence $\%{v_k\}$ approximating %$u$ strictly in $BV(\Om;\gR^2)$ and \alert{with equibounded area.}
%%This lack of degrees of freedom in the approximation of $u$ is one of the main features allowing to conclude locality property for $\widetilde {\mathcal F}_{BV}(u,\Om)$ when $\Om\subseteq\gR^2$. Notice that although the optimal minimal lifting of $u$ is unique also if $u:\Om\rightarrow\gR^2$ with $\Om\subset\gR^n$, the locality is in general false if $n\geq3$, as observed in \cite{CM}. In the case
%
{We finally recall that the locality property of the relaxed functional $A\mapsto \wid\cA_{BV}(u,A)$ fails to hold in case $n=N=3$, compare  \cite{CM}.}\footnote%
{ {Consider the vortex map $u_V(x)=x/|x|$ defined from the $3D$ unit ball $B^3$ into $\gR^3$. Roughly speaking, there are two qualitatively different Cartesian currents $T_1$ and $T_2$ that naturally arise in the relaxation procedure. They take the form $T_i=G_{u_V}+S_i$, where $S_1$ is the completely vertical current corresponding to integration of 3-forms on $\{0\}\times D^3$, where the origin $0$ is the point singularity of the vortex map and $D^3$ is the unit ball in the target space. Instead, $S_2$ corresponds to integration of 3-forms on $L\times \Sph^2$, where $L$ is a line segment connecting the origin to a point at the boundary of the domain $B^3$, whereas $\Sph^2$ is the 2-sphere in the target space given by the boundary of $D^3$.
This is the idea already contained in the counterexample from \cite{AcDM}. The main issue is that in case $n=N=3$, for $i=1,2$ we can find a sequence of smooth maps $v_k^{(i)}:B^3\to\gR^3$ converging to $u_V$ strictly in $BV$ and such that the area of the graph of $v_k^{(i)}$ converges to the mass of the current $T_i$ as $k\to\infty$, for $i=1,2$.
Notice that in case $n=N=2$, the weak convergence of the graphs to a current of the type $G_{u_V}+L\tim\Sph^1$ cannot be obtained if one imposes the strict convergence in $BV$, due to the energy concentration on the segment $L$ in the unit disk $B^2$.}
}

\smallskip

The aim of this paper is to provide some expressions for $\wid\cF_{BV}(u,\Om)$, at least in the cases when the representation formula for $\wid\cA_{BV}(u,\Om)$ is known. Among these cases we can consider $0$-homogeneous maps, some particular type of piecewise Lipschitz maps and general Sobolev maps valued in $\mathbb S^1$. These kind of maps $u$ have been considered in \cite{BSS,BSS24,Ca24,CM} where the authors proved explicit representation formulas for $\widetilde{\mathcal A}_{BV}(u,\Om)$. This is the starting point to obtain representation formulas for $\wid\cF_{BV}(u,\Om)$. Indeed we prove the following main result, where without loss of generality we choose $\Om=B^n$ the unit ball in $\gR^n$.

\bt\label{teo:intro}
Let $n\geq 2$ and let $u\in BV(B^n;\gR^2)$ be such that $\wid\cA_{BV}(u,B^n)<\i$.
Assume {that} there exists a Cartesian current $\ol T$ with underlying map $u$ and a sequence {$\{v_k\}\subset C^1(B^n;\gR^2)$} approximating $u$ strictly in $BV(B^n;\gR^2)$  such that
{
\begin{equation}
\label{hyp}
\wid\cA_{BV}(u,B^n)=\GM(\ol T)\quad{\textrm{and}}\quad G_{v_k}\rightharpoonup \overline T \quad{\textrm{as currents.}}
\end{equation}
}
Then for every continuous integrand $f$ as above
we have
$$ \wid\cF_{BV}(u,B^n)=\cF_f(\ol T)\,, $$
where $\mathcal F_f$ denotes the parametric variational integrand associated to $f$.
\et

For the precise definition of parametric variational integrand we refer to Section \ref{Sec:pvi}. We observe that the existence of a sequence {$\{v_k\}$} as in the hypothesis {\eqref{hyp}}
implies that the graph currents $G_{v_k}$ are approaching $\overline T$ in mass, namely
\begin{equation}
\label{eq:masses}
\GM(G_{v_k})\rightarrow \GM(\overline T).
\end{equation}
Then in order to prove Theorem \ref{teo:intro} we need two ingredients: The first one 
{is} continuity result, showing that if \eqref{eq:masses} holds then also
$$\mathcal E_f(v_k,B^n)\rightarrow \mathcal F_f(\overline T).$$
This is obtained {as a consequence of a result by Reshetnyak \cite{Res}, see} Proposition \ref{Pcont} below. With this into account we readily obtain that $\wid\cF_{BV}(u,B^n)\leq \cF_f(\ol T)$. To show the opposite inequality the second ingredient comes into play, namely the uniqueness of the limiting current $\overline T$ as proved in \cite{Mu22}, which entails that for any sequence {$\{u_k\}\subset C^1(B^n;\gR^2)$} approaching strictly $u$ in $BV(B^n;\gR^2)$ {and whose graph have equibounded area, up to a subsequence} there holds $G_{u_k}\wc \overline T$. So we infer, by {lower semicontinuity of $\mathcal F_f$,} that $\cF_f(\ol T)\leq \liminf_{ k}\mathcal E_f(u_k,B^n)$. We then conclude by arbitrariness of $u_k$.

Notice that in order to apply Theorem \ref{teo:intro} we need to check {that for a certain map $u$ the hypothesis \eqref{hyp} holds.}
% $\wid\cA_{BV}(u,B^n)=\GM(\ol T)$.
This requires to find out the explicit expression of $\overline T$ and to compare it with the integral representation of $\wid\cA_{BV}(u,B^n)$. In the present paper this analysis is carried on for three {kinds} of maps:
\begin{itemize}
\item[(1)] Sobolev $\mathbb S^1$-valued maps;
\item[(2)] $0$-homogeneous maps of bounded variation, {in case $n=2$};
\item[(3)] a special subclass of $BV$ consisting of piecewise Lipschitz maps, {in case $n=2$,} see Definition \ref{def:pL}.
\end{itemize}

Thanks to the results in \cite{BSS} for $n=2$ and \cite{CM} for $n\geq3$, the case (1) is well understood and it is always possible to show {that property \eqref{hyp} holds}
%equality $\wid\cA_{BV}(u,B^n)=\GM(\ol T)$
(see Section \ref{sec:sobolev}).

In the case (2), the expression of $\wid\cA_{BV}(u,{B^2})$ is related with a planar singular Plateau problem,
%namely the problem
{consisting in} finding the disk type surface of minimal area spanning a certain Lipschitz curve in $\gR^2$
{that may have}
%having also
self-intersections (see Section \ref{sec:4}).
Then to check that
%$\wid\cA_{BV}(u,B^n)=\GM(\ol T)$
{\eqref{hyp} is satisfied} requires a fine analysis of the Plateau problem, for which we give some sufficient conditions
%for this equality
{to hold} (see Corollary \ref{cor}).

Concerning point (3), piecewise Lipschitz maps are functions defined on a finite partition of the domain $B^2$, whose interfaces form a good network of curves, in the sense of Definition \ref{def:net}. For this class of {maps} the expression of $\wid\cA_{BV}(u,{B^2})$ has been studied in \cite{BSS24}, and comparing this with the structure of the current $\overline T$ {it} is possible to check {that property \eqref{hyp} holds again}
%$\wid\cA_{BV}(u,B^n)=\GM(\ol T)$ again
under suitable conditions on the behavior of $u$ around the junction points of the network.

Finally, we discuss some connection and possible application of our results to the analysis of the relaxation in $L^1$, at least in some special cases. For example, we give some sufficient {conditions} in order to characterize, for the map $u_V(x)=\frac{x}{|x|}$, the quantity  ${\widetilde{\mathcal F}_{L^1}(u_V;B^n_R)}$ when the radius $R$ of the ball $B^n_R=B_R(0)\subset\gR^n$ is big enough. 

\smallskip

The paper is organized as follows: In Section \ref{sec:notation} we introduce the notation and setting of the problem. We review the notion of relaxation of the area functional, of {Cartesian currents,} and we recall the notation in multi-vectors analysis to define the parametric polyconvex extension of the functional $\mathcal E_f(u,B^n)$ (see Section \ref{Sec:poly}). In Section \ref{sec:3} we recall the classical Reshetnyak continuity theorem and extend it in order to prove Theorem \ref{teo:intro}, which is here reformulated in Proposition \ref{Pcont} and Theorem \ref{T-cont}. In Section \ref{sec:4} we review the results from \cite{CS} about the singular Plateau problem and prove some properties of it, up to giving  the aforementioned sufficient condition in order that the mass of the current {associated} to a solution of the Plateau problem coincides with the area of the minimal {disk-type surface}. Section \ref{sec:5} is devoted to the analysis of  ${\widetilde{\mathcal F}_{L^1}(u_V;B^n_R)}$, where $u_V$ is the vortex map. In Section \ref{sec:6} {we discuss} how to apply {our} main result to obtain representation formulas for the maps in the cases (1), (2), and (3) above. We conclude with Section \ref{sec:open}, where we collect some open questions and further directions to investigate.

\section{Notation and preliminary results}\label{sec:notation}
We recall some notation on $BV$ {functions} and the relaxed area functional. We then introduce some {tools on the theory of} currents, referring to  \cite{GMS98-I,GMS98-II} for a more complete discussion.

For the basic notation, we use the symbol {$\mathcal L^n$} to denote the Lebesgue measure in {$\gR^n$;}  the symbol {$B^n_R$} denotes the open ball of radius $R>0$, centered at the origin of {$\gR^n$.} When $R=1$ we simply write {$B^n:=B_1^n$.} We also denote by $\|u\|_{p,A}$ the $L^p$-norm of {a function $u$ defined on an open set $A$. Finally, $\mathcal H^k$ denotes the $k$-dimensional Hausdorff measure.}

In what follows, $n,N\geq 2$ will be two natural numbers, and we work with currents defined in the cylinder $U=B^n\tim\RN$.
%, where $B^n$ is the unit ball in $\gr^n$.
The starting example is the $n$-current $G_u$ associated to the naturally oriented graph $\cG_u$ in $B^n\tim\RN$ of a smooth function $u:B^n\to\RN$.

\subsection{Relaxed area functional}\label{sec:21}
If $u\in C^1(B^n;\RN)$, the area functional is given by
$$ A(u,B^n):=
\int_{B^n}|M(\nabla u)|\,\dd x\,,
$$
where for any real valued $N\tim n$ matrix $G\in \gR^{N\tim n}$, we denote
\begin{equation}
\label{MG}
|M(G)|^2:= \sum_{k=0}^{\min\{n,N\}}|M_k(G)|^2 .
\end{equation}
In this formula, we let $|M_0(G)|=1$, $|M_1(G)|=|G|$, and for $2\leq k\leq \min\{n,N\}$, the squared term $|M_k(G)|^2$ is equal to the sum of the square of the determinants of all minors of order $k$ of the matrix $G$.

If e.g. $n=N=2$, then $|M_2(G)|=|\det G|$ and we have
$$ A(u,B^2)=
\int_{B^2}\{1+|\nabla u|^2+(\det\nabla u)^2\}^{1/2}\,\dd x\,.
$$
\subsubsection{BV maps }
We recall that the space of vector valued functions of bounded variations $BV(B^n;\RN)$ is the vector space of maps  $u\in L^1(B^n;\RN)$ whose distributional derivative $Du$ is a $\gR^{N\tim n}$-valued Borel measure {in $B^n$} with finite total variation, $|Du|(B^n)<\i$.
In that case, $u$ is approximately differentiable $\cL^n$-a.e. on $B^n$,
%where $\mathcal L^n$ is the Lebesgue measure in $\gR^n$,
the approximate gradient $\nabla u$ is a summable function in $L^1(B^n;\gR^{N\tim n})$, and the mutually singular decomposition holds:
$$ Du=\nabla u\,\cL^n+D^su\,,\quad |Du|(B^n)=\int_{B^n}|\nabla u|\,\dd x+|D^su|(B^n)\,. $$

In turn, the singular part $D^su$ of $Du$ splits in two mutually singular measures, $D^su=D^Cu+D^Ju$.
%, the Cantor and jump parts respectively.
{The Cantor part satisfies $|D^Cu|(B)=0$ for every Borel set $B\sb B^n$ such that $\cH^{n-1}(B)<\infty$.}
The jump part $D^Ju$ is concentrated on a {countably} $\mathcal H^{n-1}$-rectifiable set $J_u$, the jump set of $u$, and writes as
$$D^Ju=(u^+-u^-)\otimes \nu\;\dd \mathcal H^{n-1}$$
where $\nu(x)$ is {a normal unit vector} to $J_u$ at $x$, and  $u^\pm(x)$ {denote} the one-sided limits w.r.t.  $\nu(x)$ at  $x\in J_u$ (see \cite{GMS98-I}).

A sequence $\{u_h\}\sb BV(B^n;\RN)$ is said to converge weakly-$^*$ in $BV$ to a function $u\in BV(B^n;\RN)$ if $u_h\to u$ strongly in $L^1(B^n;\RN)$ and $Du_h\wc Du$
weakly-$^*$ as measures, i.e., $\langle Du_h,\vf\rangle\to \lan Du,\vf\ran$ as $h\to \i$ for every $\vf\in {C^0_b}(B^n;\gR^{N\tim n})$.
We recall that if a sequence $\{u_h\}\sb BV(B^n;\RN)$ converges strongly in $L^1$ to some function $u\in L^1(B^n;\RN)$ and $\sup_h|Du_h|(B^n)<\i$, then $u\in BV(B^n;\RN)$ and a subsequence of $\{u_h\}$ converges to $u$ weakly-$^*$ in $BV$. In that case, we also have $|Du|(B^n)\leq\liminf_h|Du_h|(B^n)$.

If $\{u_h\}\sb BV(B^n;\RN)$ converges weakly-$^*$ in $BV$ to $u\in BV(B^n;\RN)$, and
$|Du_h|(B^n)\to|Du|(B^n)$, then the sequence $\{u_h\}$ is said to converge to $u$
{\em strictly in $BV$.}

We introduce the quantity
$$ \cV(u,B^n):=\int_{B^n}\sqrt{1+|\nabla u|^2}\,\dd x +|D^s u|(B^n)\,,\quad u\in BV(B^n;\RN)\,. $$
After extending $\cV(\cdot,B^n)$ to $L^1(B^n;\gR^N)$ by setting $\cV(\cdot,B^n)=+\infty$ on $L^1(B^n;\gR^N)\setminus BV(B^n;\gR^N)$,  the functional $\cV(u,B^n)$ becomes lower-semicontinuous with respect to the $L^1(B^n;\gR^N)$ convergence. This means that  $\cV(u,B^n)\leq\liminf_h \cV(u_h,B^n)$ whenever $u_h\rightarrow u$ in $L^1(B^n;\gR^N)$ (see \cite{GMS98-I}).
\subsubsection{Relaxation of the area functional}
For a given function $u\in L^1(B^n;\RN)$, we denote by $\wid\cA_{L^1}(u,B^n)$ the relaxed area functional w.r.t. to the $L^1$-convergence, i.e.,
$$  \wid\cA_{L^1}(u,B^n)=\inf\left\{\liminf_{h\to\i}A(u_h,B^n) \mid \{u_h\}\sb C^1(B^n;\RN)\,,\,\,u_h\to u\,\,\text{in}\,\,L^1(B^n;\RN)\right\} .
$$
If $\wid\cA_{L^1}(u,B^n)<\i$, it turns out that $u\in BV(B^n;\RN)$ and
$$
\wid\cA_{L^1}(u,B^n)\geq A(u,B^n)+|D^su|(B^n)=\int_{B^n}|M(\nabla u)|\,\dd x\,+|D^su|(B^n)\geq \mathcal V(u,B^n),
$$
where $A(u,B^n)$ is the area functional, defined as above but in terms of the approximate gradient $\nabla u$.
%Notice that if $\wid\cA_{L^1}(u,B^n)<\i$ we also have $\cV(u,B^n)\leq\wid\cA_{L^1}(u,B^n)$.

In case $\wid\cA_{L^1}(u,B^n)<\i$, we also denote by $\wid\cA_{BV}(u,B^n)$ the relaxed area functional w.r.t. to the strict convergence in $BV$, i.e.,
$$  \wid\cA_{BV}(u,B^n)=\inf\left\{\liminf_{h\to\i}A(u_h,B^n) \mid \{u_h\}\sb C^1(B^n;\RN)\,,\,\,u_h\to u\,\,\text{strictly in}\,\,BV(B^n;\RN)\right\} ,
$$
so that clearly $\wid\cA_{L^1}(u,B^n)\leq \wid\cA_{BV}(u,B^n)$. Usually, this inequality might be strict and there are maps $u\in BV(B^n;\gR^N)$ such that $\wid\cA_{L^1}(u,B^n)<\infty$ and $\wid\cA_{BV}(u,B^n)=\infty$ (see \cite{BSS24}).
\subsection{Currents}  Let $m$ be a positive integer.
Given an open set $U\subset \gR^m$, for any integer $k\in \{0,\dots,m\}$ we denote by $\mathcal D^k(U)$  the
space of smooth $k$-forms with compact support in $U$ and by $\mathcal D_k(U)$ its dual, the space of $k$-dimensional currents. If $T\in \mathcal D_k(U)$ we denote by
$$\gM(T)=\sup\{T(\o) \mid \o\in\cD^k(U),\,\,\Vert\o\Vert\leq 1\}
$$
the mass of $T$. Moreover, the boundary $\partial T\in \mathcal D_{k-1}(U)$ of $T$
is the current defined by
$$\partial T(\y):=T(\dd \y)\qquad  \forall\,\omega\in \mathcal D^{k-1}(U),$$
where $\dd\y$ denotes the external differential of $\y$.
If $k=0$ by convention it holds $\partial T=0$.

Weak convergence $T_h\wc T$ of currents in $\cD_k(U)$ is defined by duality, i.e., requiring that
$$
\lim_{h\to \i} T_h(\o)=T(\o)\quad\fa\,\o\in\cD^k(U)\,.
$$
Whenever $\Phi:U\rightarrow V$ is Lipschitz and proper,
$V\subset\gR^M$ is open, and $M\geq k$, the push-forward of the current $T\in \mathcal D_k(U)$ through $\Phi$ is given by
$$ \Phi_\sharp T(\wid\o)=T(\Phi^\sharp \wid\o)\,,\quad \fa\,\wid\o\in\cD^k(V)\,,
$$
where $\Phi^\sharp \wid\o\in\cD^k(U)$ denotes the pull-back $k$-form. Therefore, $\Phi_\sharp T\in \mathcal D_k(V)$.

We say that a current $T\in \mathcal D_k(U)$ is {\it rectifiable} if there exist a {countably} $\mathcal H^k$-rectifiable set
%\footnote{$S$ is said $\mathcal H^k$-rectifiable if there are (at most) countably many Lipschitz maps $\phi_h:\R^k\rightarrow \R^n$ such that
%	$$S\subseteq N\cup\bigcup_{h=0}^{+\infty} \phi_h(\R^k),\qquad \mathcal H^k(N)=0.$$
%}
$S$ {in $U$,}
%where $\mathcal H^k$ is the $k$-dimensional Hausdorff measure,
a simple unit $k$-vector $\x(x)$ {in $\gR^m$ orienting the approximate tangent $k$-space to $S$ at $x$}
for $\mathcal H^k$-a.e. $x\in S$, and a $\cH^k\pri S$-{summable} and nonnegative map $\theta:S\rightarrow \gR$ such that
\begin{equation}\label{eq:S_rectif}
	T(\omega)=\int_S\theta(x)\langle\omega(x),\x(x)\rangle~
	\dd\mathcal H^{k}(x),\qquad \forall\,\omega\in \mathcal D^k(U),
\end{equation}
where we have denoted by $\langle\cdot,\cdot\rangle$ the standard duality pairing between $k$-covectors and $k$-vectors.
In that case, we write $T=\qu{S,\t,\x}$, {and we have $\GM(T)=\int_S {\theta}\,\dd\Ha^k<\infty$.}

A rectifiable current  $T\in \mathcal D_k(U)$ is said to be an {\it integer multiplicity current} (i.m. current) if $\theta$ takes integer values.
%, $\tau$ is tangent to $S$, and $T$ has finite mass.
{We} denote by $\cR_k(U)$ the corresponding class. If $T\in\cR_k(U)$ and in addition $\gM(\partial T)<\infty$, the boundary rectifiability theorem implies that $\pa T\in\cR_{k-1}(U)$, and $T$ is called an {\it integral current}, say $T\in\cI_k(U)$.
%
%If $k=n$ and  $S=E$ is a subset of $\gR^n$, the symbol $\qu{E}$ stands for the  integration over $E$ defined as the rectifiable $n$-current with $\theta=1$ and $\tau=e_1\wedge\dots\wedge e_n$, namely
%$$	\qu{E}(\omega)=\int_E\langle\omega(x),e_1\wedge\dots\wedge e_n\rangle~
%\dd x,\qquad \forall\omega\in \mathcal D^n(\gR^n). $$

Federer-Fleming Compactness theorem asserts that a
sequence of integral currents $T_h\in \mathcal I_k(U)$ with $\sup_h (\gM(T_h)+
\gM(\partial T_h))<\infty$ admits a
subsequence weakly converging to an integral current $T \in
\cI_k(U)$.

If $S$ is a smooth oriented $k$-manifold embedded in $U$, the symbol $\qu{S}$ denotes the current naturally associated by integration of $k$-forms on $S$ in the sense of Differential Geometry. Therefore, if $S$ has a naturally oriented smooth boundary $\pa S$, by Stokes theorem we have
$$ \pa\qu{S}(\y)=\qu S(\dd\y)=\int_S \dd\y=\int_{\pa S}\y\,, \quad\fa\,\y\in\cD^{k-1}(U)\,.
$$
\subsubsection{Multi-vectors}
We let $\LL_n\gR^{n+N}$ denote the space of $n$-vectors $\x$ in
$\gR^{n+N}\simeq\Rn\tim\RN$.
Symbols $\{e_1,\dots, e_n\}\subset\LL_1\gR^n$ and $\{\varepsilon_1,\dots,\varepsilon_N\}\subset\LL_1\gR^N$ denote the canonical bases of $1$-vectors in $\gR^n$ and $\RN$, and
$\{\dd x^1,\dots,\dd x^n\}$ and $\{\dd y^1,\dots,\dd y^N\}$ the dual bases of $1$-covectors  in $\gR^n$ and $\RN$, respectively.

We denote by $\x^{\ol 0 0}$ the coefficient of $e_1\wedge\cdots\wedge e_n$ in the expression of a $n$-vector $\x\in\LL_n\gR^{n+N}$, that is $\x^{\ol 0 0}=\lan \dd x^1\wedge\cdots\wedge \dd x^n,\x\ran$.
We also denote by $\SS$ the class of simple $n$-vectors
in $\Lambda_n\gR^{n+N}$, and we let
%that is the ones generated by the wedge product of $n$ 1-vectors. We also let
$$
\ba{rcl}
\LL_1 & := & \{\x\in \LL_n\gR^{n+N}\mid \x^{\ol 00}=1\}\,, \\
\SS_1 & := & \{\x\in \SS\mid \x^{\ol 00}=1\}\,, \\
\SS_+ & := & \{\x\in \SS\mid \x^{\ol 00}>0\}\,. \ea
$$

For any real valued $N\tim n$ matrix $G\in \gR^{N\tim n}$, we let
$$ \vect M(G):=(e_1+G e_1)\wdg\cdots\wdg(e_n+G e_n)\in\LL\nolm_n\RnN , $$
where
$$
G e_i =\sum_{j=1}^N\,G^j_i\,\e_j\,,\quad i=1,\ldots,n\,.
$$
Then $\vect M(G)\in \SS_1$, and the unit simple $n$-vector
$$\x_G:=\frac{\vect M(G)}{| \vect M(G)|} $$
identifies the plane graph of $G$ in $\gR^{n+N}$, and in fact orients such a
$n$-plane. We have $\x_G\in \SS_+$, with $\x_G^{\ol 0 0}=| \vect M(G)|^{-1}$ and $| \vect M(G)|=|M(G)|$, see \eqref{MG}.

Notice moreover that for every $\x\in \Sigma_1$ there exists a unique matrix $G_\x$ in $
\gR^{N\tim n}$ such that
$$\x=\vect M(G_\x)\,. $$
{We finally recall that $\LL_1$ is the convex envelop of $\Sigma_1$.}
\begin{definition}\label{def:poly}
A function $f:\gR^{N\tim n}\to\gR\cup\{+\i\}$ is said to be {\em polyconvex} if there exists a convex function $g: {\LL_1} \to\gR\cup\{+\i\} $ such that
$$
f(G)=g(\vect M(G))\quad\fa\, G\in \gR^{N\tim n}
$$
or equivalently
$$
g(\x)=f(G_\x)\quad\fa\, \xi\in \Sigma_1\,.
$$
\end{definition}
\subsubsection{Graph currents}
The $n$-current $G_u$ carried by the graph of a smooth map $u:B^n\to\RN$ is defined by the integration of compactly supported and smooth $n$-forms $\o$ in $\cD^n(B^n\tim\RN)$ on the naturally oriented graph $\cG_u$ of $u$. By the area formula, we have
\beq\label{Gu}
{G_u(\o)}=\int_{B^n}\lan \o(x,u(x)),\vect M(\nabla u(x)) \rangle\,\dd x
\eeq
and hence the mass of $G_u$ is equal to the area functional, i.e.,
\beq\label{areamass}
\GM(G_u)=A(u,B^n)\,.
\eeq
Furthermore, by Stokes theorem for every $\y\in\cD^{n-1}(B^n\tim\RN)$
we have 
$$ {\pa G_u(\y) := G_u(\dd\y)}=\int_{\cG_u}\dd\y=\int_{\pa \cG_u} \y =0\,, $$
and hence $\pa G_u=0$ on $\cD^{n-1}(B^n\tim\RN)$. Therefore, if $|M(\nabla u)|\in L^1(B^n)$, it turns out that $G_u$ is an integral $n$-current with null boundary in $B^n\tim\RN$.
\smallskip

We introduce the class $\cA^1(B^n;\RN)$ consisting of all functions $u\in L^1(B^n;\RN)$ which are approximately differentiable $\cL^n$-a.e. and such that $|M(\nabla u)|\in L^1(B^n)$, where $\nabla u$ is the approximate gradient.

If $u\in \cA^1(B^n;\RN)$, the graph current $G_u$ is well-defined in an approximate sense by formula \eqref{Gu}, and it is again an i.m. rectifiable current in $\R_n(B^n\tim\Rn)$, with finite mass satisfying equation \eqref{areamass}. However, in general the null-boundary condition $\pa G_u=0$ is violated, and we may have $\GM(\pa G_u)=\i$.
\subsubsection{Cartesian currents}
Cartesian currents naturally arise by applying Federer-Fleming Compactness theorem to sequences of graph currents of smooth maps. More precisely, if $\{u_h\}\sb C^1(B^n;\RN)$ satisfies
$$
\sup_h\left(\GM(G_{u_h})+\Vert u_h\Vert_{\i,B^n}  \right)<\i \,,
$$
recalling that $\pa G_{u_h}=0$ for every $h$, we find a not relabeled subsequence and an i.m. rectifiable current $T\in\cR_n(B^n\tim\RN)$ such that $G_{u_h}\wc T$ weakly as currents. Moreover, by lower semicontinuity we have
$$ \GM(T)\leq\liminf_{h\to\i}\GM(G_{u_h})<\ii $$
and, by stability of the null-boundary condition, we have $\pa T=0$. Furthermore, there exists a map of bounded variation $u_T$ in $\cA^1(B^n;\RN)$ such that $u_h\to u_T$ in $L^1(B^n;\RN)$. Therefore, we can decompose
\beq\label{decT}
T=G_{u_T}+S_T\,,\quad \GM(T)=\GM(G_{u_T})+\GM(S_T)
\eeq
where $S_T$ is a {\em vertical} current in $\cR_n(B^n\tim\RN)$, in the sense that
$$
{ S_T(\vf(x,y)\,\dd x^1\wedge\cdots\wedge\dd x^n)}=0\quad \fa\,\vf\in C^\i_c(B^n\tim\RN)\,.
$$
Finally, the support of $T$ is contained in the compact set 
$$\cK_R:= \overline B^n\times \overline B^N_R$$ given by the closure (in $B^n\times \gR^N$) of $B^n\tim B^N_R$ for some radius $R>0$, that is, $T(\o)=0$ for every $\o\in\cD^n(B^n\tim\RN)$ that is null in $\cK_R$.

The class of Cartesian currents with compact support is identified by the previous properties. Therefore, in this paper we say that $T\in
\cart(B^n\tim\RN)$ if:
\ben
\item $T\in\cR_n(B^n\tim\RN)$;
\item $T$ has compact support contained in $\cK_R$ for some $R>0$;
\item $T$ satisfies the null boundary condition $\pa T=0$;
\item $T$ can be decomposed as in \eqref{decT} for some
function $u_T\in\cA^1(B^n;\RN)$ and some vertical current $S_T$ in $\cR_n(B^n\tim\RN)$;
\item $u_T$ is a function of bounded variation in $BV(B^n;\RN)$, with $\Vert u_T\Vert_{\infty,B^n}\leq R$.
\een

It turns out that the class $\cart(B^n\tim\RN)$ is closed along weakly converging sequences with equibounded mass and support contained in some $\cK_R$.
\subsection{Relaxed area and currents}
 For a function $u\in L^1(B^n;\RN)$ such that $\wid\cA_{L^1}(u,B^n)<\i$ we introduce the class
$$\cT_u:=\{T\in \cart(B^n\tim\RN):T=G_u+S_T,\;\text{for some vertical current }S_T\in \mathcal R_n(B^n\times \gR^N)\},$$
i.e., the family
of Cartesian currents
 with underlying map $u$.
We recall that if $u\in L^1(B^n;\RN)$ is such that $\wid\cA_{L^1}(u,B^n)<\i$, then $u\in\cA^1(B^n;\RN)$ and the class $\cT_u$  is non-empty. Moreover we have
$$ \wid\cA_{L^1}(u,B^n)=\inf\{\wid\cA(T)\mid T\in\cT_u\}\,, $$
where for any $T\in\cT_u$ we have let
$$  \wid\cA(T):=\inf\left\{\liminf_{h\to\i}A(u_h,B^n) \mid \{u_h\}\sb C^1(B^n;\RN)\,,\,\,G_{u_h}\wc T\,\,\text{weakly in}\,\,\cD_n(B^n\tim\RN)\right\}\,.
$$

Assume now that $\wid\cA_{BV}(u,B^n)<\i$. Then, arguing as in \cite{Mu22}, it turns out that there exist currents $T$ in $\cT_u$ whose action on forms with exactly one differential in the vertical direction only depends on $u$ and $Du$. More precisely, we denote by $\cT_u^{BV}$ the elements in $\cT_u$ such that for every $i\in\{1,\ldots,n\}$, $j\in\{1,\ldots,N\}$, and $\f\in C^\i_c(B^n\tim\RN)$ we have
$$ { T(\f(x,y)\wih{\dd x^i}\wedge\dd y^j)}=(-1)^{i-1}\int_{B^n}\left(\int_0^1\f(x,u^\t(x))\,\dd\t\right)\dd (Du)^j_i \,, $$
where $u^\t(x)=\t\,u^+(x)+(1-\t)\,u^-(x)$ if $x\in J_u$, and $u^\t(x)$ agrees with a precise representative $\cH^{n-1}$-a.e. on $B^n\sm J_u$.

In fact, if $\wid\cA_{BV}(u,B^n)<\i$ we have
$$ \wid\cA_{BV}(u,B^n)=\inf\{\wid\cA(T)\mid T\in\cT_u^{BV}\}\,.$$
However, in general the class $\cT_u^{BV}$ contains a unique element only if $N=2$, for every $n\geq 2$, see \cite{Mu22}.

By lower semicontinuity, we thus have
$$ \wid\cA_{L^1}(u,B^n)\geq \inf\{\gM(T)\mid T\in\cT_u\}\,,\quad \wid\cA_{BV}(u,B^n)\geq \inf\{\gM(T)\mid T\in\cT_u^{BV}\}\,,$$
but the strict inequality holds true in presence of topological obstructions that are not seen by the homology. This is the case of the double-eight example (see Remark \ref{rem_28}),
given by the 0-homogeneous extension $u_8$ of the map $\vf_8\in \Lip(\Sph^1;\gR^2)$
\begin{equation}
\label{ottodoppio}
\vf_8(\cos\t,\sin\t):=\left\{
\ba{ll}
(-1+\cos 4\t,\sin 4\t)  & \If\quad 0\leq\t<\p/2       \\
(1-\cos 4\t,\sin 4\t)   & \If\quad \p/2\leq\t<\p      \\
(-1+\cos 4\t,-\sin 4\t) & \If\quad \p\leq\t<3\p/2   \\
(1-\cos 4\t,-\sin 4\t)  & \If\quad 3\p/2\leq\t<2\p  \\
\ea\right.
\end{equation}
whose image covers an ``eight" figure twice but with opposite orientation.
{The Sobolev function $u_8$ belongs to the class $\A^1(B^2;\gR^2)$ and its graph current satisfies $\pa G_{u_8} =0$, whence $G_{u_8}\in\cart(B^2\tim\gR^2)$, but (see, e.g., \cite[Section 6]{Sc})}
$$ \wid\cA_{L^1}(u_8,B^2)>\gM(G_{u_8})\,,\quad \wid\cA_{BV}(u_8,B^2) > \gM(G_{u_8})\,,\quad \gM(G_{u_8})=\mathcal V(u_8,B^2)\,. $$
\subsection{Polyconvex energies}\label{Sec:poly}
In this paper we consider nonnegative and continuous integrands
$$f:B^n\tim\RN\tim
\gR^{N\tim n}\to\gR_+ , $$
where $\gR_+:=[0,+\i)$, satisfying the following
properties:
\ben\item[$(a)$] the function $G\mapsto f(x,u,G)$ is polyconvex in $
\gR^{N\tim n}$ for every $(x,u)\in B^n \tim\RN$;\footnote%
{This means that the function $G\mapsto f(x,u,G)$ can be written as a convex function $g(M(G))$ of the minors of $G\in\gR^{N\tim n}$, see Definition \ref{def:poly}.
}
\item[$(b)$]
{there exists a real constant $C>0$ such that
$$ 0\leq f(x,u,G)\leq C\,|M(G)| $$}
for every $(x,u,G)\in B^n \tim\RN\tim
\gR^{N\tim n}$.
\een
The corresponding energy on smooth functions $u:B^n\to\RN$ is given by
\begin{equation}\label{Fu}
\cE_f(u,B^n):=\int_{B^n}f(x,u,\nabla u)\,\dd x\,.
\end{equation}

Consider the map
$$\ul f: B^n\tim\RN\tim \SS_1\to\ol\gR_+ $$
defined by
$$\ul f(x,u,\x):=f(x,u,G_\x)\,, $$
where, for every $\x\in \Sigma_1$, we have denoted by $G_\x$ the unique matrix in $
\gR^{N\tim n}$ such that $\x=\vect M(G_\x)$.

Taking $x,u$ as parameters, we let
$$\ol  f_{x,u}:\Lambda_n\gR^{n+N}\to [0,+\i]
$$
be given by
$$\ol f_{x,u}(\x):=\left\{\ba{ll}\renewcommand{\arraystretch}{1.8}
\x^{\ol 00}\ul f({x,u},\x/\x^{\ol 00})
%= \x^{\ol 00} f_{x,u}(G_\x)
&
\quad{\mx{\text if}}\quad \x\in\SS_+ \\
+\i & \quad{\mx{\text otherwise}}\,. \ea\right. $$
\begin{definition}
The parametric polyconvex l.s.c. envelope of the integrand $G\mapsto
f_{x,u}(G):=f(x,u,G)$ is given by the convex l.s.c. envelope of
$\x\mapsto \ol f_{x,u}(\x)$, namely
$$
F_{x,u}(\cdot):=\Gamma C \ol f_{x,u}(\cdot)\,.
$$
\end{definition}
In fact, by the continuity of the integrand $f$, it turns out that the function
$(x,u,\x)\mapsto F_{x,u}(\x)$ is lower semicontinuous in all variables and
convex in $\x$ for every $(x,u)\in B^n \tim\RN$. Moreover, we have
$$ F_{x,u}(\x)=F_f(x,u,\x) $$
for every $(x,u,\x)\in B^N\tim\RN\tim\Lambda_n\gR^{n+N}$, where
$$
\ba{r}\renewcommand{\arraystretch}{1.8}
F_f(x,u,\x):=\sup\{\phi(\xi)  \mid  \phi:\Lambda_n\gR^{n+N}\to\ol\gR_+\,,\,\phi\,\,{\mx{\text linear}},\qquad \\
\phi(\vect M(G))\leq f(x,u,G)\quad\fa\, G\in \gR^{N\tim n} \} \,.
\ea
$$
\subsubsection{Parametric variational integrals}
\label{Sec:pvi}
For our purposes, if $T\in\cD_{n}( B^n \tim\RN)$ has finite mass, we denote by $\Vert T\Vert$ the mass norm, i.e., the finite Borel measure {given by}
$$\Vert T\Vert(V):=\GM(T\pri V) $$
for each {open} set $V\sb B^n \tim\RN$, by $T=\Vert T\Vert\pri\vect T $ the corresponding Radon-Nikodym decomposition,
and by $\pi$ and $\wih\pi$ the orthogonal projections onto the two factors in $\Rn\tim\RN$.
%, see \rif{enorm}. Here \,$T$\, is identified with the
%$\gR^{1+Nn}$-valued linear functional
%$$T:=\bigl(T^{\ol 00},(T^{\ol ij})_{\gR^{Nn}}\bigr)\,,\qquad i=1,\ldots
%n\,,\quad j=1,\ldots N \,. $$
% where for every \,$\phi\in C^\i_0( B^n \tim\calY)$
% $$ T^{\ol 00}(\phi):=T(\phi\,dx)\,,\qquad
% T^{\ol ij}(\phi):=T(\phi\,\wih{dx^i}\wdg dy^j)\,. $$
%
\begin{definition}
The {\em parametric variational integral} associated to the
integrand $f$ is defined by
\beq\label{7pF}
\cF_f(T,B\tim\RN):=\int_{B\tim\RN}F_f\bigl(\p(z),\wih\p(z),\vect
T(z)\bigr)\,\dd\Vert T\Vert(z)
\eeq
for every Borel set $B\sb B^n $,
%where $F_f(x,u,\xi)$ is ...
%given by $\rif{7F}$,
and we let
$$ \cF_f(T):=\cF_f(T, B^n\tim\RN)\,. $$
\end{definition}

{Notice that for currents $T$ in $\cart(B^n\tim\RN)$, if $f(G)=|M(G)|$ we have $\cF_f(T)=\gM(T)$. Moreover, by property (b)
we infer that
$$
\cF_f(T)\leq C\,\GM(T). $$
%for every current $T\in\cD_{n}(B^n\tim\RN)$ with finite mass.
%
Finally,} by the construction it turns out that the following lower
semicontinuity property holds (see \cite{GMS98-II}).
\bp\label{Plsc} Let
$\{T_h\}\sb\cD_{n}(B^n\tim\RN)$ be such that $\sup_h\GM(T_h)<\i$ and $T_h\wc T$ weakly in $\cD_{n}$.
Then
$$ \cF_f(T)\leq\liminf_{h\to\i}\cF_f(T_h). $$
\ep

If $T\in\cR_n(B^n\tim\RN)$, and $T=\qu{S,\t,\x}$, compare \eqref{eq:S_rectif}, we have $\vect T=\x$ and $\Vert T\Vert=\t\,\cH^n\pri S$, so that
$$ \cF_f(T)= \int_{S}\t(z)\,F_f\bigl(\p(z),\wih\p(z),\x(z)\bigr)\,\dd\cH^n(z)\,. $$

If e.g. $T=G_u$ for some $u\in\cA^1(B^n;\RN)$, using that $\t=1$ and $S=\cG_u$, the {\em rectifiable graph} of $u$, recalling that
$$
\x(x,y)=\frac{\vect M(\nabla u(x))}{ |\vect M(\nabla u(x))| } \,,$$
for $\cH^n$-almost every $(x,y)\in\cG_u$, and that $F_f\bigl(x,u,\vect M(G) \bigr)=f(x,u,G)$, by the area formula we obtain
$$ \cF_f(G_u)=
\int_{\cG_u} F_f\bigl(x,u(x),\vect M(\nabla u(x)) \bigr)\,\frac 1{|\vect M(\nabla u(x) |}\dd\cH^n(z) =\int_{B^n}f(x,u(x),\nabla u(x))\,\dd x\,, $$
so that extending the notation in \eqref{Fu} we have
$$ \cF_f(G_u)=\cE_f(u,B^n)=\int_{B^n}f(x,u,\nabla u)\,\dd x\,\quad\fa\, u\in\cA^1(B^n;\RN)\,, $$
where $\nabla u$ is the approximate gradient of $u$.
More generally, if $T\in\cart(B^n\tim\RN)$, and $T=G_u+S_T$, by the mass norm decomposition $\Vert T\Vert =\Vert G_u \Vert + \Vert S_T \Vert$ we have
$$ \cF_f(T)=\cE_f(u,B^n)+ \cF_f(S_T)\,. $$
%
%
%
%%%%%
%If $T\in\cR_n(B^n\tim\RN)$ is of the type $T=\qu{S,\t,\tau}$, compare \eqref{eq:S_rectif}, we denote
%$$
%\cF_f(T):=\int_{B^n\tim\RN} \t(x,u)\,F_f(x,u,\tau(x,u))\,\dd \cH^n\pri S \,.
%$$
%Notice that $\cF_f(T)=+\i$ if there exists a subset of points in $S$ with positive $\cH^n$ measure where the multiplicity $\t$ is positive but $\tau^{\ol 0 0}<0$. This is not the case of Cartesian currents.
%Moreover, we point out that property (a) implies that for every $u\in\cA^1(B^n;\RN)$ we have
%$$ \cF_f(G_u)= \int_{B^n}f(x,u,\nabla u)\,\dd x =:\cE_f(u,B^n)\,.
%$$
%
%Moreover, if $\{T_h\},T_\i\sb \cR_n(B^n\tim\RN)$ and $T_h\wc T_\i$ weakly as currents, then
%$$\cF_f(T_\i)\leq\liminf_{h\to\i}\cF_f(T_h)\,. $$

Let now $\{u_h\}\sb C^1(B^n;\RN)$ be a sequence of smooth functions with equibounded energies and uniform norms,
\begin{equation}
\label{bound}
\sup_h\left( \cE_f(u_h,B^n) + \Vert u_h \Vert_{\i,B^n}\right)<\infty\,.
\end{equation}
In order to apply Federer-Fleming Compactness theorem, the energy density needs to satisfy the following regularity property:
\ben
\item[$(c)$] there exists a real constant $c>0$ such that
$$ f(x,u,G)\geq c\,|M(G)| $$
for every $(x,u,G)\in B^n \tim\RN\tim \gR^{N\tim n}$.
\een

In that case, in fact,
{for every sequence $\{u_h\}\sb C^1(B^n;\RN)$ satisfying \eqref{bound},}
we obtain the mass bound $\sup_h\GM(G_{u_h})<\infty$, so that
we find a not relabeled subsequence and a current $T\in\cart(B^n\tim\RN)$ such that
$G_{u_h}\wc T$. Moreover, the support of $T$ is contained in $\cK_R$, where  $R=\sup_h\Vert u_h \Vert_{\i,B^n}$.
\br
Assume now that $N=2$ and that we know in addition that the sequence $\{u_h\}$ converges strictly in $BV$ to some $u\in BV(B^n;\RN)$. Then, instead of
{the regularity} property (c), in order to apply Federer-Fleming Compactness theorem and to conclude as above, it suffices to require the
{weaker lower bound:}
\ben
\item[$(c')$]
{there exists a real constant $c>0$ such that
$$ f(x,u,G)\geq c\,|M_2(G)| $$
for every $(x,u,G)\in B^n \tim\gR^2\tim\gR^{2\tim n}$.}
\een
In that case, moreover, the underlying map $u_T$ in \eqref{decT} is equal to $u$.
%\er
%
%\begin{remark}
%	{\rm

\smallskip

A well studied functional satisfying condition $(c')$ and not $(c)$ is the {\it total variation of the Jacobian} for maps $u:B^2\rightarrow \gR^2$. Precisely, one considers
$$
f(x,u,G)=|\det G|,\qquad {(x,u,G)\in B^2\tim\gR^2\tim\gR^{2\tim 2}.}
$$
		The corresponding energy, for smooth $u:B^2\rightarrow \gR^2$, is given by
$$	%	\begin{equation}
			TVJ(u,B^2):=\int_{B^2}|\det(\nabla u)|\,\dd x.
$$ %		\end{equation}
%	}
See \cite{DP} and references therein for relaxation of this functional in the Sobolev setting. 
\er
%%%%%%%
%
\section{An extension of Reshetnyak continuity theorem}\label{sec:3}
The following continuity theorem is due to
Reshetnyak \cite{Res}, compare Thm.~1 in Sec.~1.3.4 of
\cite{GMS98-II}. Here, for $m\geq2$, {and for} any $\gR^m$-valued Radon measure
$\m$ defined on an open set $U\sb\gR^{n+N}$, we denote
by $\vect\m$ its Radon Nikodym derivative with respect to the
total variation $|\m|$, and by $\m_k\wc\m$ the weak-$^*$
convergence in the sense of the measures.
\bt\label{TR2}{\bf(Reshetnyak).} Let $G(z,p)$ be a
non-negative continuous function defined in $U\tim\gR^m$ {and}
satisfying the following properties:
\ben
\item $G(z,\cdot)$ is
positively homogeneous of degree one for every $z$;
\item $G(\cdot,p)$ is uniformly bounded for every $p\in\Sph^{m-1}$;
\item $G(z,\cdot)$ is {\em essentially convex} for every $z$, i.e.,
$$ G(z,p+q)\leq G(z,p)+G(z,q)\qquad\fa\,p,q\in\gR^m, $$
where the equality holds if and only if $q=\l p$ for some
$\l\geq 0$.
\een
Let $F(z,p)$ be a non-negative continuous function that is
{positively} homogeneous of degree one in $p$ for every $z$ and that
satisfies
$$ 0\leq F(z,p)\leq c_1\,G(z,p)+c_2\qquad\fa\,(z,p)\in U\tim\gR^m$$
for some absolute constants $c_1,c_2>0$. Then we have
$$ \lim_{k\to\i}\int_U F(z,\vect\m_k(z))\,d|\m_k|=\int_U F(z,\vect\m(z))\,d|\m| $$
provided that $\m_k$, $\m$ are $\gR^m$-valued Radon
measures on $U$ satisfying
$$ \m_k\wc\m\,,\qquad \int_U G(z,\vect\m_k(z))\,d|\m_k|\to\int_U G(z,\vect\m(z))\,d|\m|\qquad {\text{as }}\,k\to\i\,.$$
\et
\subsection{A new continuity result}
The following result extends to the $BV$ setting a property observed by Acerbi--Dal Maso \cite{AcDM}, who gave a sufficient condition to the strong convergence in $W^{1,1}$.

\bt\label{Tstrict}
Let $\{u_h\}\sb BV(B^n;\RN)$ be a sequence weakly-$^*$ converging in $BV$ to some $u\in BV(B^n;\RN)$, and assume that $\cV(u_h,B^n)\to \cV(u,B^n)$ as $h\to \i$. Then we have
\begin{equation}\label{conv_thesis}
	\lim_{h\to \i}\left( \int_{B^n}|\nabla u_h-\nabla u|\,\dd x +|D^s u_h|(B^n)
	\right) =|D^s u|(B^n)
\end{equation}
and hence $\{u_h\}$ converges to $u$ strictly in $BV$.
\et
\bpf For a given $v\in BV(B^n;\RN)$, we let $\m(Dv)$ denote the $\gR^{1+nN}$-valued measure in $B^n$ given by
$$ \m(Dv):=\left(\calL^n,\,(Dv)_i^j\right)\,,\quad i=1,\ldots,n\,,\quad j=1,\ldots,N\,. $$
By the mutual singularity of the absolutely continuous and singular components of {$Dv$,} we have
$$ \m(Dv)=\Bigl(1,\,\nabla_iv^j\Bigr)\,\calL^n+\Bigl(0,\,(D^sv)_i^j\Bigr)\,,\quad
|\m(Dv)|=\sqrt{1+|\nabla v|^2}\,\cL^n+|D^sv|\,,
$$
so that $|\m(Dv)|(B^n)=\cV(v,B^n)$. Therefore, by the assumptions we have that $\m(D{u_h})\wc \m(Du)$ weakly-$^*$ as measures and $|\m(D{u_h})|(B^n)\to |\m(Du)|(B^n)$ as $h\to \i$.
Given $\vf\in C^0_c(B^n;\gR^{nN})$, define the continuous function $\psi:B^n\tim{\gR^{1+nN}}\to\gR$ by $\psi(x,\z)=|\hat\z -\z_0\,\vf(x)|$, if $\z=(\z_0,\hat \z)\in \gR\tim\gR^{nN}\simeq \gR^{1+nN}$.
Arguing as in \cite[Lemma 6.3]{AcDM}, by Theorem \ref{TR2} we infer that
$$ \lim_{h\to \i}\left( \int_{B^n}|\nabla u_h-\vf|\,\dd x +|D^s u_h|(B^n)
\right) = \int_{B^n}|\nabla u-\vf|\,\dd x +|D^s u|(B^n)\,. $$
By a density argument, the latter property holds true for any $\vf\in L^1(B^n;\gR^{nN})$, and hence we can choose $\vf=\nabla u$. The last property in the statement follows from the second triangle inequality; indeed, on the one hand {by lower semicontinuity}
$$\liminf_{h\to\i}\int_{B^n}(|\nabla u_h|-|\nabla u|)\dd x+|D^su_h|(B^n)-|D^su|(B^n)\geq 0.$$
On the other hand, by \eqref{conv_thesis} and the triangle inequality
$$\limsup_{h\to\i}\!\int\limits_{B^n}\!(|\nabla u_h|-|\nabla u|)\dd x+|D^su_h|(B^n)-|D^su|(B^n)
%=\limsup_{h\to\i}\int_{B^n}(|\nabla u_h|-|\nabla u|-|\nabla u_h-\nabla u|)\dd x\leq 0,
{
\leq\lim_{h\to\i}\!\int\limits_{B^n}\! |\nabla u_h-\nabla u|\dd x+|D^su_h|(B^n)-|D^su|(B^n)= 0
}
$$
that {give} the thesis.
\epf
\par A first consequence of Theorem \ref{Tstrict} is the following property inspired by a result from \cite{Sc}.

\bt\label{teo3.3} Let $u\in L^1(B^n;\RN)$ be such that $\wid\cA_{L^1}(u,B^n)<\i$, and assume that
$\wid\cA_{L^1}(u,B^n)=\cV(u,B^n)$. Then we also have $\wid\cA_{BV}(u,B^n)=\cV(u,B^n)$.
\et
\bpf Let $\{u_h\}\sb C^1(B^n;\RN)$ be such that $u_h\to u$ strongly in $L^1$ and $A(u_h,B^n)\to\wid\cA_{L^1}(u,B^n)$. Then, possibly passing to a not relabeled subsequence we have $u_h\wc u$ weakly-$^*$ in $BV$. By the lower semicontinuity of the functional $u\mapsto\cV(u,B^n)$ we have
$$ \cV(u,B^n)=\wid\cA_{L^1}(u,B^n)=\lim_{h\to\i}A(u_h,B^n)\geq \limsup_{h\to\i}\cV(u_h,B^n)
\geq \liminf_{h\to\i}\cV(u_h,B^n)\geq \cV(u,B^n)\,, $$
so that we obtain $\cV(u_h,B^n)\to \cV(u,B^n)$ as $h\to \i$. Therefore, Theorem \ref{Tstrict} implies that $u_h\to u$ strictly in $BV$.
{We thus obtain}
$$ \wid\cA_{BV}(u,B^n)\leq\lim_{h\to\i}A(u_h,B^n)=\cV(u,B^n)=\wid\cA_{L^1}(u,B^n)\leq \wid\cA_{BV}(u,B^n)$$
and the assertion readily follows.
\epf
\par As a consequence, for every $u\in BV(B^n;\RN)$ we have
\begin{equation}\label{gaps}
\wid\cA_{BV}(u,B^n)>\cV(u,B^n) \Longrightarrow \wid\cA_{L^1}(u,B^n)>\cV(u,B^n) \,.
\end{equation}

\begin{definition}
	{Whenever $\wid\cA_{L^1}(u,B^n)>\int_{B^n}|M(\nabla u)|\dd x=\gM(G_u)$ we say that there is an {\em energy gap} for the relaxed  area functional in $L^1$. Similarly, there is an energy gap for the relaxed area functional in $BV$ if $\wid\cA_{BV}(u,B^n)>\gM(G_u)$.}
\end{definition}

\subsection{Further remarks}
Assume now that $u\in BV(B^n;\RN)$ is such that
\begin{equation}\label{1dim}
|M(\nabla u)|=\sqrt{1+|\nabla u|^2} ,
%\qquad \qquad |M_2(\nabla u)|=0\,,
\end{equation}
i.e., all the subdeterminants of order $k\geq2$ of the matrix $\nabla u$ are null.
This happens if e.g. $u$ takes values in a 1-dimensional set of $\RN$, by the area formula.

We thus have $u\in\cA^1(B^n;\RN)$ and
$$ A(u,B^n)=\cV(u,B^n)\,. $$

{Therefore, if $\wid\cA_{BV}(u,B^n)<\i$
and if \eqref{1dim} takes place,
by \eqref{gaps} it turns out that} the occurrence of an energy gap for the relaxed area in $BV$ implies an energy gap for the relaxed area in $L^1$, as observed in \cite{Sc} for the double eight example.

Notice that Theorem \ref{teo3.3} implies that if there is no energy gap for either $\wid\cA_{L^1}(u,B^n)$ or $\wid\cA_{BV}(u,B^n)$, then
$$\wid\cA_{L^1}(u,B^n)=\wid\cA_{BV}(u,B^n).$$
The converse is false in general, i.e., if $\wid\cA_{L^1}(u,B^n)=\wid\cA_{BV}(u,B^n)$ then energy gaps might occur. Consider for example the vortex map  given by  $u_V(x)=x/|x|$ where $n=N=2$. After replacing $B^2$ by a bigger ball ${B=B^2_R}\subset\gR^2$ (or equivalently, rescaling the map $u_V$), the results in \cite{BSSbis} imply that, if $R$ is big enough, then
$$\wid\cA_{L^1}(u_V,B)=\mathcal V(u_V,B)+\pi,$$
where the energy gap $\pi$ coincides with the area enclosed by the unit circle, in turn being equal to
$$
{
P\Bigl(\frac{x}{|x|}\Bigr) = \inf\left\{ \int_{B^2}|\det \nabla v|\,\dd x \,:\,
v\in \text{{\rm Lip}}(B^2;\gR^2),\;v=\frac{x}{|x|} \text{ on }\partial B^2 \right\}\,.
}
$$
From this and \cite{BSS} (see also Theorem \ref{T-BSS} below), we infer
$$\wid\cA_{L^1}(u_V,B)=\wid\cA_{BV}(u_V,B).$$
\subsection{A result on Cartesian currents}
We now briefly discuss an extension of Theorem \ref{Tstrict} to Cartesian currents. We show that the mass convergence is equivalent to an a priori stronger energy convergence.
\bp Let $\{T_h\}\sb\cart(B^n\tim\RN)$ be such that $\sup_h\GM(T_h)<\i$ and
the support of every $T_h$ is contained in $\cK_R$ for some $R>0$. Then, there exists a not relabeled subsequence and a current $T_\i\in\cart(B^n\tim\RN)$ such that
$T_h\wc T_\infty$ weakly in $\cD_n(B^n\tim\RN)$, $T_\infty$ is supported in $\mathcal K_R$, and
\beq\label{lsc}
\GM(T_\i)\leq\liminf_{h\to\i}\GM(T_h)\,.
\eeq
Moreover, writing $T_h=G_{u_h}+S_{T_h}$ for every $h\in{\mathbb N}\cup\{\i\}$, the mass convergence $\gM(T_h)\to \gM(T_\i)$ as $h\to\i$ is equivalent to the energy convergence
\beq\label{M-strong}
\lim_{h\to \i}\left( \int_{B^n}|M(\nabla u_h)-M(\nabla u_\i)|\,\dd x +\GM(S_{T_h}) \right) =\GM(S_{T_\i})\,.
\eeq
\ep
\bpf The first assertion follows directly from Federer-Fleming Compactness theorem.
As to the equivalence property, the first implication follows by arguing as in the proof of Theorem \ref{Tstrict} above. Conversely, if \eqref{M-strong} holds, recalling that
$$
\GM(T_h)=\int_{B^n}|M(\nabla u_h)|\,\dd x +\GM(S_{T_h})
$$ for every $h\in\ol\Nat$, by the second triangle inequality we have
$$ \limsup_{h\to\i} \GM(T_h)\leq \GM(T_\i)\,, $$
and hence the mass convergence follows from the lower semicontinuity inequality \eqref{lsc}.
\epf

Notice that in the previous proposition, if $T_h=G_{u_h}$ for some smooth function $u_h$, %for every $h$ and $T_\i\in \cT^{BV}_u$ for some $u$ with $\wid\cA_{BV}(u,B^2)<\i$, by %the convergence $G_{u_h}\wc T_i$ with $A(u_h,B^n)\to\gM(T_\i)$
we infer that
$$ \lim_{h\to \i}\int_{B^n}|M(\nabla u_h)-M(\nabla u_\i)|\,\dd x =\GM(S_{T_\i})\,. $$
However, we cannot conclude that $u_h\to u$ strictly in $BV$, even if we knew that ${T_\i}\in \cT_{u_\i}^{BV}$.

Therefore, the latter result does not seem to give new information concerning the relaxed area functionals, {when compared to} the one discussed in the previous section.
\subsection{A continuity property for energies on currents}
In the following, we let $f$ be a continuous integrand
$$f:B^n\tim\RN\tim
\gR^{N\tim n}\to\gR_+ $$
satisfying properties $(a)$ and $(b)$ from Sec. \ref{Sec:poly}. We then let $F_f$, $\cE_f$, and $\cF_f$ be defined as before.

For our purposes, we now recall how the following continuity property is obtained from Theorem \ref{TR2}:
\bp\label{Pcont} {
Let $\{T_h\}\sb \cart(B^n\tim\RN)$ be such that
%$\sup_h\gM(T_h)<\i$ and
$T_h\wc T$ weakly in $\cD_n(B^n\tim\RN)$ to
some $T\in \cart(B^n\tim\RN)$.} If
$\gM(T_h)\to\gM(T)$, then $\cF_f(T_h)\to\cF_f(T)$.
\ep
\bpf We set $U=B^n\tim\RN$,
$z=(x,u)$, and $m=c(n,N)$ where $c(n,N)$ is the dimension of the ordered vector corresponding to ${\vect M}(G)$ for some $G\in \gR^{N\tim n}$.
Set now
$$G(z,p):=|p|\,,\qquad F(z,p):=F_f(x,u,p) $$
if $u\in\RN$ and $p$ is identified with the components
of an $n$-vector $\xi\in\Lambda_n(\Rn\tim \RN)$ with  $\x^{\ol 00}\geq 0$.
By suitably extending $G$ and $F$, it is readily checked
that we can apply Theorem~\ref{TR2}. Since the convergence of the corresponding measures $\m_{(T_h)}\wc \m_{(T)}$ reduces to the weak convergence
$T_h\wc T$ in $\cD_n(B^n\times \gR^N)$, whereas
$$ \int_U G(z,\vect\m_{(T)}(z))\,\dd |\m_{(T)}|=\gM(T)\,,\qquad \int_U F(z,\vect\m_{(T)}(z))\,\dd |\m_{(T)}|=\cF_f(T)\,, $$
the assertion readily follows. \epf
\subsection{Relaxed energies}
For a given function $u\in L^1(B^n;\RN)$, we denote by $\wid\cF_{L^1}(u,B^n)$ the relaxed functional w.r.t. to the $L^1$-convergence, i.e.,
$$  \wid\cF_{L^1}(u,B^n)=\inf\left\{\liminf_{h\to\i}\cE_f(u_h,B^n) \mid \{u_h\}\sb C^1(B^n;\RN)\,,\,\,u_h\to u\,\,\text{in}\,\,L^1(B^n;\RN)\right\}.
$$
If $f$ also satisfies property (c) from Sec. \ref{Sec:poly}, then $\wid\cF_{L^1}(u,B^n)<\i$ if and only if $\wid\cA_{L^1}(u,B^n)<\i$. In that case, moreover, we have
$$ \wid\cF_{L^1}(u,B^n)=\inf\{\wid\cF(T)\mid T\in\cT_u\}, $$
where for any $T\in\cT_u$ we have let
$$  \wid\cF(T):=\inf\left\{\liminf_{h\to\i}\cE_f(u_h,B^n) \mid \{u_h\}\sb C^1(B^n;\RN)\,,\,\,G_{u_h}\wc T\,\,\text{weakly in}\,\,\cD_n(B^n\tim\RN)\right\}.
$$
In addition, if $\wid\cA_{L^1}(u,B^n)<\i$ we also denote
$$  \wid\cF_{BV}(u,B^n)=\inf\left\{\liminf_{h\to\i}\cE_f(u_h,B^n) \mid \{u_h\}\sb C^1(B^n;\RN)\,,\,\,u_h\to u\,\,\text{strictly in}\,\,BV(B^n;\RN)\right\},
$$
so that clearly $\wid\cF_{L^1}(u,B^n)\leq \wid\cF_{BV}(u,B^n)$, and $\wid\cF_{BV}(u,B^n)<\i$ if and only if $\wid\cA_{BV}(u,B^n)<\i$, where
$$ \wid\cF_{BV}(u,B^n)=\inf\{\wid\cF(T)\mid T\in\cT_u^{BV}\}. $$
\par
The following representation result holds true for the relaxed energy in the strict convergence, in the case of codimension $N=2$. In that case, we have seen that we can require the weaker lower bound $(c')$ instead of $(c)$ from Sec. \ref{Sec:poly}.
\bt\label{T-cont} Let $n\geq 2$ and let $u\in BV(B^n;\gR^2)$ be such that $\wid\cA_{BV}(u,B^n)<\i$.
Assume that there exists a Cartesian current $\ol T$ in $\cT_u$ such that $\wid\cA_{BV}(u,B^n)=\GM(\ol T)$. Then for every continuous integrand $f$ satisfying properties
{$(a),\,(b),\,(c')$ in Sec. \ref{Sec:poly}}
we have
$$ \wid\cF_{BV}(u,B^n)=\cF_f(\ol T)\,. $$
\et
\bpf We recall from \cite{Mu22} that if $N=2$ and $\wid\cA_{BV}(u,B^n)<\i$, there exists a unique Cartesian current $T_u$ in the class $\cT_u^{BV}$, so that by the strict convergence we have $\overline T=T_u$.
Since $\wid\cF_{BV}(u,B^n)<\i$, by closure-compactness and lower semicontinuity, see Proposition \ref{Plsc}, we infer that
$$
\wid\cF_{BV}(u,B^n)\geq\inf\{\cF_f(T)\mid T\in{\cT_u^{BV}}\}=\cF_f(\ol T)\,.
$$
Moreover, there exists a sequence $\{u_h\}\sb C^1(B^n;\gR^2)$ such that $u_h\to u$ strictly in $BV$ and $A(u_h,B^n)\to\gM(\ol T)$. Possibly passing to a not relabeled subsequence, we infer that $G_{u_h}\wc \ol T$ and $\gM(G_{u_h})\to\gM(\ol T)$.
By the continuity result in Proposition \ref{Pcont}, applied with $T_h=G_{u_h}$ and $T=\ol T$, we infer that $\cE_f(u_h,B^n)=\cF_f(T_h)\to\cF_f(\ol T)$ as $h\to\i$ and the proof is complete.
\epf
%
%%%%%%%%%
%
\section{Generalized curves and a planar Plateau problem}\label{sec:4}
Let
$\gamma\in BV([a,b];\gR^2)$ be a given function of bounded variation defined on a closed nondegenerate interval $[a,b]$. We denote by ${D\g}$ the distributional derivative of $\g$.

If $L_\gamma:=|\dot \gamma|([a,b])$ we introduce the map
\begin{equation}\label{s_gamma}
	s_\gamma(t)=\frac{1}{L_\gamma+(b-a)}\big(|{D\g}|([a,t))+(t-a)\big), \qquad \qquad \forall t\in [a,b],
\end{equation}
which will be strictly increasing with {discontinuities in the jump set $J_\gamma$ of $\gamma$;} we observe that $s_\gamma$ is right continuous and
$$s_\gamma(t_1)-s_\gamma(t_2)\geq \frac{t_1-t_2}{L_\gamma+(b-a)},\qquad \qquad {a\leq t_2 < t_1\leq b,}$$
and so it follows that, denoting $ t_\gamma:=s_\gamma^{-1}:[0,1]\rightarrow [a,b]$  the {continuous and non decreasing} inverse of $s_\gamma$ that is constant on $[s_\gamma(t^-),s_\gamma(t^+)]$ for all $t\in S_\gamma$, we have
\begin{equation}
	t_\gamma(s_1)-t_\gamma(s_2)=|t_\gamma(s_1)-t_\gamma(s_2)|\leq (s_1-s_2)(L_\gamma+(b-a)),\qquad \qquad 0\leq s_2\leq s_1\leq1.
\end{equation}
This shows that $t_\gamma$ is Lipschitz continuous with Lipschitz constant $L_\gamma+(b-a)$.
\begin{definition}\label{defrep}
	Given $\gamma\in BV([a,b];\gR^2)$ we define $\overline \gamma:[0,1]\rightarrow \gR^2$ as
	\begin{equation}\label{gencurve}
		\overline \gamma(s)=\begin{cases}
			\displaystyle
			\frac{\gamma(t)^+(s-s_\gamma(t)^-)+\gamma(t)^-(s_\gamma(t)^+-s)}{s_\gamma(t)^+-s_\gamma(t)^-}&\text{if }s\in [s_\gamma(t)^-,s_\gamma(t)^+],\\
			\gamma(t_\gamma(s))&\text{otherwise.}
		\end{cases}
	\end{equation}
	We call $\overline \gamma$ the {\em generalized curve} associated to $\gamma$.
\end{definition}

According to \cite[Proposition 3.6]{Sc}, the curve $\overline \gamma$ is Lipschitz continuous with Lipschitz constant less than or equal to $L_\gamma+(b-a)$.

By \cite[Proposition 3.6]{Sc} (see also \cite{Ca24}) the following property holds: if $\gamma_k\in BV([a,b];\gR^2)$ are  converging strictly to $\gamma\in BV([a,b];\gR^2)$ then
$\overline \gamma_k$ tend to $\overline \gamma$ uniformly. In the next section we will use this fact to investigate a planar Plateau problem $\overline P(\gamma)$ associated to any $BV$ function $\gamma$.
	Given a  Lipschitz curve $\varphi:\mathbb S^1\rightarrow \gR^2$ we define
$$P(\varphi):=\inf\left\{\int_{B^2}|Jv|dx:v\in \text{{\rm Lip}}(B^2;\gR^2),\;v=\varphi \text{ on }\partial B^2\right\}.$$
This gives the area of a disk type solution of a two dimensional Plateau problem spanning the curve $\varphi(\mathbb S^1)$.
\subsection{A planar Plateau problem}
The relaxation of $P$ on $BV(\mathbb S^1;\gR^2)$ is given by, for any $\gamma\in BV(\mathbb S^1;\gR^2)$,
\begin{equation}\label{plateau_rel}
\overline P(\gamma):=\inf \{\liminf_{k\rightarrow \infty}P(\varphi_k):\varphi_k\in \text{{\rm Lip}}(B^2;\gR^2),\; \varphi_k\rightarrow \gamma \text{ strictly in }BV \}.
\end{equation}
Up to setting $a=0$, $b=2\pi$ and considering {functions} with
{$\gamma(0)=\gamma(2\pi)$,} we can associate to each $\gamma\in BV(\mathbb S^1;\gR^2)$ its reparametrization {$\overline \gamma$} as in Definition \ref{defrep}.
Since $P(\varphi)$, for $\varphi\in \text{{\rm Lip}}(\mathbb S^1;\gR^2)$, is invariant under reparametrization of $\varphi$, we deduce that {$P(\gamma)=P(\overline \gamma)$,} and thanks to the fact that $P$ is continuous under uniform convergence in $\text{{\rm Lip}}(\mathbb S^1;\gR^2)$ (see \cite{Ca24}),  by the property above we infer
\begin{equation}\label{plateau_rel=plateau}
	\overline P(\gamma)=\lim_{k\rightarrow \infty}P(\overline \varphi_k)=P(\overline \gamma).
\end{equation}
This shows that the computation of $\overline P$ can always be led back to $P$.

{For a given $\varphi\in \text{{\rm Lip}}(\mathbb S^1;\gR^2)$, the} quantity $P(\varphi)$ can be deduced from the geometry of $\varphi(\mathbb S^1)$ and the way $\varphi$ parameterizes it, as shown in \cite{CS}. To describe how $P(\varphi)$ is computed, we make the assumption that
\begin{itemize}
	\item[(P0)] The open set $\gR^2\setminus \varphi(\mathbb S^1)$ has a finite number of connected components.
\end{itemize}

Under this hypothesis, we select a family of points $\{q_1,q_2,\dots,q_K\}$, each  taken in a given bounded connected component $U_i$, $i=1,\dots,K$, and another  fixed point $q_0\in \varphi(\mathbb S^1)$ called base point. We then consider the homotopy group $\pi_1(\gR^2\setminus \{q_1,q_2,\dots,q_K\})$ of loops based at $q_0$, and a basis of it $\{\sigma_1,\sigma_2,\dots,\sigma_K\}$, where $\sigma_i$ is a simple Lipschitz  loop based at $q_0$ winding (with counter-clockwise orientation) one time the point $q_i$ and zero times the other points $q_j$, $j\neq i$. We denote by $\sigma_i^{-1}$ the same loop as $\sigma_i$ but with opposite orientation (in such a case $\sigma_i$ and $\sigma_i^{-1}$ are said inverse to each other). Then it is well-known that  $\pi_1(\gR^2\setminus \{q_1,q_2,\dots,q_K\})$ is isomorphic to the free group $\mathbb F(K)$ on $K$ letters, where the elements of it can be identified with words (or chains) containing letters in the alphabet $\{\sigma_1,\sigma_2,\dots,\sigma_K,\sigma_1^{-1},\sigma_2^{-1},\dots,\sigma_K^{-1}\}$.

A word in $\mathbb F(K)$ has two {\it extremal letters} (the first and the last one), and we say that  a word in $\mathbb F(K)$ is  {\it reduced} if it does not contain $\sigma_i$ and $\sigma_i^{-1}$ as adjacent letters (for some $i$). For example, $\sigma_1\sigma_2\sigma_3^{-1}$ and $\sigma_2\sigma_1^{-1}\sigma_2\sigma_2$ are reduced words, whereas $\sigma_1\sigma_2\sigma_2^{-1}$ and $\sigma_2\sigma_1\sigma_1^{-1}\sigma_3$ are not. Given a word $\beta$, we can modify it into a reduced word $\beta'$ by erasing adjacent couples $\sigma_i\sigma_i^{-1}$. We say that two words $\beta$ and $\gamma$ are {\it equivalent}, and we write $\beta\sim\gamma$, if they can be modified into the same reduced word. Finally, two words $\beta$ and $\gamma$ are said to be {\it conjugate} if there is a word $\sigma$ such that $\beta\sim\sigma\gamma\sigma^{-1}$.
For example, the words $\sigma_2\sigma_1\sigma_1^{-1}\sigma_3$ and $\sigma_2\sigma_2^{-1}\sigma_2\sigma_3$ can be reduced to $\sigma_2\sigma_3$, and so $\sigma_2\sigma_1\sigma_1^{-1}\sigma_3\sim\sigma_2\sigma_2^{-1}\sigma_2\sigma_3$; also, $\sigma_1^{-1}\sigma_2\sigma_1$ and $\sigma_2$ are conjugate to each other. The {\it conjugate class} of a word $\gamma$ is the family containing all the words conjugate to $\gamma$ and is denoted by $[\gamma]\subset \mathbb F(K).$ We denote by $0$ the empty word; any word belonging to $[0]$ is said a {\it null word}.
\medskip

To every curve $\varphi$ satisfying hypothesis (P0) we can associate in a unique way a conjugate class in $\mathbb F(K)$, denoted by $[\varphi]$. Indeed $\varphi$ is an element of $\pi_1(\gR^2\setminus \{q_1,q_2,\dots,q_K\})$ and for this reason can be represented by a word in $\mathbb F(K)$. In general there is more than one word representing $\varphi$, but all these representations are conjugate to each other (see \cite{CS}).

We now describe an algebraic operation that can be acted on words and that is needed to describe how to compute $P(\varphi)$.

\begin{definition}
	For $i\in \{1,\dots,K\}$, an $i$-monoid
	is a word conjugate to either $\sigma_i$ or $\sigma_i^{-1}$.
	Let $\gamma\in \mathbb F(K)$ and let $(k_1 ,\dots, k_K)\in \mathbb N^K$ be a $K$-tuple of natural
	numbers. We say that $\gamma'\in \mathbb F(K)$ is a $(k_1 ,\dots, k_K)$- injection in $\gamma$ if the word $\gamma'$ is in the conjugate class of a word obtained by inserting, for all $i = 1,\dots,K$, $k_i$ times an $i$-monoid into the word $\gamma$.
\end{definition}		

For instance, given $\gamma=\sigma_1 \sigma_2\in \mathbb F(3)$ then $\sigma_1 \sigma_3 \sigma_2$ is a $(0, 0, 1)$-injection in $\gamma$. The
word $\sigma_1$ instead is a $(0, 1, 0)$-injection  in $\gamma$ (we inserted $\sigma_2^{-1}$ into the word of $\gamma$).
But it is also a $(0, 3, 0)$-injection if we suitably insert $\sigma_2 , \sigma_2^{-1} , \sigma_2^{-1}$ in $\gamma$. So there is not a unique way for seeing a word $\gamma'$ as an injection in $\gamma$. The null word is a
$(1, 1, 0)$-injection in $\gamma$ (as well as a $(1 + 2h, 1 + 2h, 2h)$-injection for all $h \in \mathbb N$).
Furthermore,  $\sigma_1 \sigma_3 \sigma_2$ is a $(0, 0, 1)$-injection in $\gamma$, as it can be obtained by inserting
the monoid $\sigma_2^{-1}\sigma_3\sigma_2$ at the end of the word representing $\gamma$; moreover $\sigma_1\sigma_2\sigma_1^{-1}\sigma_2^{-1}\in\mathbb  F(2)$ is
a $(1, 0)$-injection in $\sigma_1$, as it is obtained from $\sigma_1$ by inserting the $1$-monoid $\sigma_2\sigma_1^{-1}\sigma_2^{-1}$.
Finally, $0$ is a $(2,0)$-injection (and also a $(0,2)$-injection) in $\sigma_1\sigma_2\sigma_1^{-1}\sigma_2^{-1}$ (and viceversa).

Observe also that if $\gamma''$ is a $(k_1 , \dots, k_K)$-injection in $\gamma'$ and $\gamma'$  is a $(h_1 , \dots, h_K)$-injection in $\gamma$,
then $\gamma''$ is a $(k_1 + h_1 ,\dots, k_K + h_K )$-injection in $\gamma$.

\begin{remark}\label{rem_1}
	We notice that, given a word $\gamma\in \mathbb F(K)$, if $\gamma$ contains (for all $i=1,\dots,K$) $n_i$ times the letter $\sigma_i$ or $\sigma_i^{-1}$, then the null word is a $(n_1,\dots,n_K)$-injection in $\gamma$. Indeed, it is sufficient to insert, for all letters in $\gamma$, the monoid consisting of the corresponding inverse letter beside each of them.
\end{remark}

\begin{definition}
{	Given $\varphi$ satisfying hypothesis (P0), let $\gamma\in [\varphi]$; we introduce the class
	$$\text{{\rm Ad}}(\gamma):=\{(k_1,\dots,k_K)\in \mathbb N^K:0 \text{ is a $(k_1,\dots,k_K)$-injection in $\gamma$}\}.$$}
\end{definition}	
By Remark \ref{rem_1} $\text{{\rm Ad}}(\gamma)$ is never empty. Moreover, it is easy to see that $\text{{\rm Ad}}(\gamma)$ does not depend on the choice of $\gamma\in [\varphi]$.
For this reason, we will write
$$\text{{\rm Ad}}(\gamma)=\text{{\rm Ad}}([\varphi]).$$

The main result in \cite{CS} reads as follows:

\begin{theorem}\label{teo_plateau}
	Let $\varphi\in \text{{\rm Lip}}(\mathbb S^1;\gR^2)$ satisfy hypothesis (P0); then
	$$P(\varphi)=\min \Big\{\sum_{i=1}^Kk_i\mathcal L^2(U_i):(k_1,\dots,k_K)\in \text{{\rm Ad}}([\varphi])\Big\},$$
	where we recall that $U_i$ is the $i$-th {bounded connected component of}
$\gR^2\setminus \varphi(\mathbb S^1)$.
\end{theorem}

\begin{definition}
Let $(k_1,\dots,k_K),(k'_1,\dots,k'_K)\in\text{{\rm Ad}}(\varphi)$ be two injections.
{If $k_i\leq k'_i$ for all $i=1,\dots,K$, we write $(k_1,\dots,k_K)\prec (k'_1,\dots,k'_K)$. Moreover,} we say that $(k_1,\dots,k_K)$ is a minimal injection if whenever $(k'_1,\dots,k'_K)\prec (k_1,\dots,k_K)$ then $(k'_1,\dots,k'_K)= (k_1,\dots,k_K)$.
{Finally,} given a word $\gamma\in \mathbb F(K)$, we denote by $\text{{\rm Ad}}_{\min}(\g)$ the class of minimal injections in $\text{{\rm Ad}}(\g)$.
\end{definition}		

We observe that  a minimal injection always exists; however it does not need to be unique. For example, suppose $\sigma_1\sigma_2\sigma_1^{-1}\sigma_2^{-1}=\gamma\in [\varphi]\subseteq \mathbb F(2)$, then $(2,0)$ and $(0,2)$ are two distinct minimal injections.
\subsection{Relationship with currents}
We now turn to the Plateau problem $P(\varphi)$ and observe that, given $\varphi\in \text{{\rm Lip}}(\mathbb S^1;\gR^2)$, by
{the Constancy theorem} there exists a unique $2$-current $S_\varphi$ in $\mathcal D_2(\gR^2)$, with compact support, such that
$$\partial S_\varphi=\varphi_\sharp\qu{\mathbb S^1}=:\jump{\varphi}.$$
It is easy to see that this current can be written as
\begin{equation}\label{currentSvarphi}
S_\varphi=\sum_{i=1}^K h_i\jump{U_i},
\qquad \qquad h_i:=\text{{\rm wind}}(\varphi;q_i),
\end{equation}
where $h_i:=\text{{\rm wind}}(\varphi;q_i)$ is the winding number of the curve $\varphi$ around the point $q_i$.

\begin{remark}\label{rem_28}
Notice that in general $P(\varphi)$ does not coincide with the mass of the current $S_\varphi$ in \eqref{currentSvarphi}, i.e., $k_i$ in Theorem \ref{teo_plateau} is not necessarily equal to $|h_i|$. Consider indeed the double-eight curve $\varphi_8$
{in \eqref{ottodoppio};} this is defined as
the concatenation
$$\varphi_8=\sigma_1^{-1}\sigma_2^{-1}\sigma_1\sigma_2,$$
where $\sigma_1$ and $\sigma_2$ are, respectively, circles of radius one centered in $q_1:=(1,0)$ and $q_2:=(-1,0)$ (hence meeting at the origin) running in counter-clockwise sense. The corresponding current $S_{\varphi_8}$ is null by \eqref{currentSvarphi}. However, according to Theorem \ref{teo_plateau}, $P(\varphi_8)=2\pi$, namely $P(\varphi_8)$ coincides with two times the area of one of the two {disks}, because the minimal injections in such a case are $\text{{\rm Ad}}_{\min}([\varphi_8])=\{(2,0),(0,2)\}$. We will next prove that if the minimal injection is unique, then $P(\varphi)=\gM(S_\varphi)$.
\end{remark}

Under the hypothesis that there is a unique minimal injection
$(k_1,\dots,k_K)\in { \text{{\rm Ad}}_{\min} } (\varphi)$, from Theorem \ref{teo_plateau} it readily follows that
$$P(\varphi)=\sum_{i=1}^Kk_i\mathcal L^2(U_i).$$

We will now show that in such a case
{we have $k_i=|h_i|$ for every $i$.}
%, up to a sign, $k_i=h_i$.
To begin with, we prove the following lemma showing that the order of making injections commutes:

\begin{lemma}
	If $\gamma\in \mathbb F(K)$ can be reduced to $\gamma'$ by inserting in it a $i$-monoid and then a $j$-monoid, then $\gamma$ can be reduced to $\gamma'$ by inserting in it a $j$-monoid and then a $i$-monoid.
\end{lemma}
\bpf
	By \cite[Remark 2.30]{CS} any monoid can be equivalently inserted at the end of the word. Hence, assume that
	$$\gamma'\sim\gamma (\alpha \sigma_i\alpha^{-1})(\beta\sigma_j\beta^{-1}),$$
	for some words $\alpha$ and $\beta$, then also
	$$\gamma'\sim\gamma (\alpha \sigma_i\alpha^{-1}\beta\sigma_j\beta^{-1}\alpha \sigma_i^{-1}\alpha^{-1}) (\alpha \sigma_i\alpha^{-1})$$
	where between parentheses any monoid has been emphasized, showing the thesis.
\epf

%
%\begin{lemma}
%	If $(k_1,\dots,k_K)\in\text{{\rm Ad}}(\varphi)$, and $\gamma\in [\varphi]$, then the null word can be obtained by inserting, for all $i=1,\dots,K$,  $k_i$ times a reduced $i$-monoid of the form $\sigma_i$ or $\sigma_i^{-1}$.
%\end{lemma}
%\bpf
%s
%\epf

\begin{theorem}\label{teo2.9}
	Assume that $\varphi\in \text{{\rm Lip}}(\mathbb S^1;\gR^2)$ satisfies hypothesis (P0) and that there exists a unique minimal injection $(k_1,\dots,k_K)\in\text{{\rm Ad}}_{\min}([\varphi])$. Then, denoting by {$\#(\sigma_i)$ and $\#(\sigma_i^{-1})$} the number of times that {$\sigma_i$ and $\sigma_i^{-1}$ appear} in $\gamma$, for all $\gamma\in [\varphi]$ it holds $|\#(\sigma_i)-\# (\sigma_i^{-1})|=k_i.$ As a consequence
	$$S_\varphi=\sum_{i=1}^K\alpha_i k_i\qu{U_i},$$
	for suitable signs $\alpha_i\in \{-1,1\}$.
\end{theorem}
\bpf Let $\gamma\in \mathbb F(K)$ be a reduced word in $[\varphi]$. We prove the first assertion in the statement.
	From this the second one will follow, since it is straightforward that $\#(\sigma_i)-\# (\sigma_i^{-1})=h_i:=\text{{\rm wind}}(\varphi;q_i)$.

	By definition of $\text{{\rm Ad}}([\varphi])$, for any $i=1,\dots,K$, there exist two natural numbers $k_i^\pm$ with $k^+_i+k^-_i=k_i$ such that we can inject $k_i^+$ times a suitable conjugate of $\sigma_i$ and $k_i^-$ times a conjugate of $\sigma_i^{-1}$ in $\gamma$ (for all $i$) and obtain the null word. In the null word there holds, for all $i$,
	$l_i:=\#(\sigma_i)=\# (\sigma_i^{-1}).$ So we deduce that in $\gamma$ it holds
	$$\#(\sigma_i)+k_i^+=l_i\qquad \qquad \#(\sigma_i^{-1})+k_i^-=l_i.$$
	We therefore deduce that $|\#(\sigma_i)-\#(\sigma_i^{-1})|=|k_i^+-k_i^-|$.
	To prove the  claim we can then show that either $k_i^+$ or $k_i^-$ is zero.
	
	Assume by contradiction that both are positive, and suppose without loss of generality that $k_i^+\geq k_i^->0$.
	In this case, we can inject in the word $\gamma$ the  inverse of a letter, beside this letter. As in Remark \ref{rem_1} (denoting by $n_j$ the number of appearance of $\sigma_j$ or $\sigma_j^{-1}$ in $\gamma$) we can do this for all letters but the ones coinciding with $\sigma_i$ and $\sigma_i^{-1}$. In this way, after reducing the obtained word, we will have the final word consisting of exactly $k_i^+-k_i^-$ times the letter $\sigma_i$. At this point it is enough to inject   $k_i^+-k_i^-$ times the letter $\sigma_i^{-1}$ to obtain the null word.
	
	This procedure shows that we can obtain, from $\gamma$, the null word by a $(n_1,\dots,(k_i^+-k_i^-),\dots,n_K)$-injection. Hence $(n_1,\dots,(k_i^+-k_i^-),\dots,n_K)\in \text{{\rm Ad}}([\varphi])$ is such that the $i$-th entry is $k_i^+-k_i^-< k_i^+<k_i$, which in turn implies the existence of minimal injections in $\text{{\rm Ad}}([\varphi])$ whose $i$-th entry is smaller than $k_i$.
{We thus}
%In particular, we
obtain a contradiction with the hypothesis that $(k_1,\dots,k_K)$ was the unique minimal injection in $\text{{\rm Ad}}([\varphi])$.	
\epf

As a consequence of the previous result, if there is a unique minimal injection in $\text{{\rm Ad}}(\varphi)$, then we have
\begin{equation}\label{equality_masscurrent}
	P(\varphi)=\gM(S_\varphi).
\end{equation}
\subsection{A characterization}
We will now characterize the words $\gamma\in\mathbb F(K)$ for which $\text{{\rm Ad}}_{\min}(\gamma)$ consists of a unique element.

\bdf\label{defH}
Let $\g\in \mathbb F(K)$. We say that $\g$ satisfies property (H) if there exists a minimal injection $(k_1,\ldots,k_N)$ in $\text{{\rm Ad}}_{\min}(\g)$ such that $|  \#(\s_i^+) - \#(\s_i^-) |=k_i$ for every $i=1,\ldots,N$, where we have denoted $\s_i^+:=\s_i$ and $\s_i^-=\s_i^{-1}$, and $\#(\sigma_i^{\pm})$ is the number of appearances of $\sigma_i^{\pm}$ in the reduced representation of $\gamma$.
\edf

%
%  for every $i$. For every $\g\in P(K)$, we let $\#(\s_i^\pm)\in\gn$ be the numbers of letters $\s_i^\pm$ appearing in the reduced representation in $[\g]$.
%For any simbol $(k_1,\ldots,k_K)\in \text{{\rm Ad}}(\g)$, writing $k_i=k_i^++k_i^-$, for every $i$ there exists $\ell_i\in\gn$ such that $ k_i^\pm+\ell_i=\#(\s_i^\pm)$, where both alternatives are assumed.
%
\bt\label{Tmin} The class $\text{{\rm Ad}}(\g)$ admits a unique minimal injection if and only if $\g$ satisfies property $(H)$.
%In that case, moreover, we have that
%$|  \#(\s_i^+) - \#(\s_i^-) |=k_i$ for every $i=1,\ldots,N$ where  $(k_1,\ldots,k_N)$ is the unique element in $\text{{\rm Ad}}_{\min}(\g)$.
\et
\bpf If $\g$ has a unique element $(k_1,\ldots,k_K)$ in $\text{{\rm Ad}}_{\min}(\g)$, then it satisfies property (H) by Theorem \ref{teo2.9}.
{Conversely,}
%Viceversa,
assume that $\g$ satisfies property (H), and let  $(k_1,\ldots,k_K),\, (\hat k_1,\ldots,\hat k_K)\in \text{{Ad}}_{\min}(\g)$ be two minimal injections, where without loss of generality we may and do assume that $|  \#(\s_i^+) - \#(\s_i^-) |=\hat k_i$ for every $i$.
Since in the null word the number of appearances of $\sigma_i^+$
{is equal to the one of $\sigma_i^-$,} for every $i$ there exist  $\ell_i,\,\hat\ell_i\in\gn$ such that both the alternatives
$$\#(\s_i^\pm)+k_i^\pm=\ell_i,\qquad \qquad \#(\s_i^\pm)+\hat k_i^\pm=\hat\ell_i$$
hold, implying that
\begin{equation}\label{diff}
	k_i^+-k_i^-=\#(\s_i^-)-\#(\s_i^+)=\hat k_i^+-\hat k_i^-,
\end{equation}
for all $i$. Hence $\hat k_i=|  \#(\s_i^+) - \#(\s_i^-) |=| k_i^+-k_i^-|\leq k_i^++k_i^-=k_i$ for all $i$, which gives a contradiction with {the} minimality of $(k_1,\ldots,k_K)$ unless $\hat k_i=k_i$.
%
%
%namely
%$$ \ell_i-k_i^\pm=\hat\ell_i-\hat k_i^\pm .$$ This yields to equations
%\beq\label{diff}
%\hat k_i^+-k_i^+=\hat\ell_i-\ell_i=\hat k_i^--k_i^-
%\eeq
%for every $i$.
%%Furthermore, we have $\hat k_i^+\cdot \hat k_i^-=0$ for every $i$.
%Assume now that
%$|  \#(\s_i^+) - \#(\s_i^-) |\neq k_i$ for some $i$. Recalling that $\#(\s_i^+),\,k_i^\pm\in\gn$, in that case we have
%$\hat k_i = |  \#(\s_i^+) - \#(\s_i^-) |< k_i$ and hence
%we obtain a contradiction to the minimality property of $(k_1,\ldots,k_N)$.
%Therefore, we have that $|  \#(\s_i^+) - \#(\s_i^-) |=k_i$ for every $i$.
%
%\par Moreover, for every $i$ at least one of the two following equivalent properties is satisfied:
%$$k_i^-=0 \iff \hat k_i^-=0\,,\quad k_i^+=0 \iff \hat k_i^+=0 \,. $$
%In the first case, we obtain $\hat k_i=\hat k_i^+=k_i^-=k_i$, and in the second one
%we obtain $\hat k_i=\hat k_i^-=k_i^-=k_i$, for every $i$, whence $(k_1,\ldots,k_N)= (\hat k_1,\ldots,\hat k_N)$.
The uniqueness of the minimal injection readily follows.
\epf

Since \eqref{equality_masscurrent} holds if there is a unique minimal injection, we readily obtain:
\begin{corollary}\label{cor}
	If $\varphi\in  \text{{\rm Lip}}(\mathbb S^1;\gR^2)$ {satisfies} (P0) and  $\gamma\in [\varphi]$ {satisfies} property (H), then \eqref{equality_masscurrent} holds.
\end{corollary}
We now deal with the homology class of the chain $\g$, that is given by
$$
\qu{\g}:= \sum_{i=1}^K \left(\#(\s_i^+) -\#(\s_i^-)\right)\,\qu{\s_i}\,.
$$
On account of property \eqref{diff}, it turns out that
$$ \#(\s_i^+) -\#(\s_i^-)=k_i^--k_i^+\quad\fa\,i\in\{1,\dots,K\} $$
independently of the choice of $(k_1,\ldots,k_K)$ in $\text{{\rm Ad}}_{\min}(\vf)$.
%, where we have used the previous notation concerning $k_i=k_i^++k_i^-$.
We thus have
\beq\label{homol}
\qu{\g}= \sum_{i=1}^K \left( k_i^--k_i^+\right)\,\qu{\s_i}\quad\fa\,(k_1,\ldots,k_K)\in \text{{\rm Ad}}_{\min}(\vf)\,.
\eeq
By property \eqref{homol} we infer that $h_i=\text{{\rm wind}}(\varphi;q_i)=k_i^--k_i^+$, i.e.,
$$
S_\vf= \sum_{i=1}^K \left( k_i^--k_i^+\right)\,\qu{U_i}
%\quad\fa\,(k_1,\ldots,k_N)\in \text{{\rm Ad}}_{\min}(\vf)
$$
that clearly implies
$$
\gM(S_\vf)=\sum_{i=1}^K | k_i^+-k_i^- |\,\cL^2(U_i)\,.$$
We are now ready to show the converse of Corollary \ref{cor}:
%%
%\subsection{Plateau's problem}
%Let now $\vf\in \Lip(\cS^1,\gr^2)$ satisfy (P0) and let $P(K)$ be the corresponding free group. Let
%$$ P(\vf)=\inf\left\{ \int_{B^2}|\det Dv|\,\dd x \mid v\in \Lip(B^2;\gr^2),\,\,v_{\vert\pa B^2}=\vf\right\}\,.
%$$
%Moreover, let $S_\vf$ be the unique current in $\cR_2(\gr^2)$ such that $\pa S_\vf=\vf_\flat\qu{\gs^1}$. We have:
%

\bt\label{Tgap} Let $\vf\in \Lip(\cS^1,\gr^2)$ satisfy $(P0)$. Then,
\beq\label{ineq}
\gM(S_\vf)\leq P(\vf)\,.
\eeq
Furthermore, the equality sign holds true in formula \eqref{ineq} if and only if the homotopy class $[\vf]$ satisfies property $(H)$, and in that case
we have
\beq\label{nogap}
\gM(S_\vf) = P(\vf)= \sum_{i=1}^K k_i\,\cL^2(U_i)
\eeq
for the unique element $(k_1,\ldots,k_K)$ in $\text{{\rm Ad}}_{\min}{(\vf)}$.
\et
\bpf In general, on account of the previous notation, and using that the class $\text{{\rm Ad}}_{\min}{(\vf)}$ has a finite number of elements, we have
$$
P(\vf)=\min\left\{\sum_{i=1}^K k_i\,\cL^2(U_i) \,:\, (k_1,\ldots,k_K)\in \text{{\rm Ad}}_{\min}{(\vf)} \right\}\,.
$$
We thus may and do choose $(k_1,\ldots,k_K)\in \text{{\rm Ad}}_{\min}(\vf)$ in such a way that
$$ P(\vf)=\sum_{i=1}^K k_i\,\cL^2(U_i)\,. $$
%Furthermore, by property \eqref{homol} we infer that
%$$
%S_\vf= \sum_{i=1}^N \left( k_i^+-k_i^-\right)\,\qu{U_i}
%%\quad\fa\,(k_1,\ldots,k_N)\in \text{{\rm Ad}}_{\min}(\vf)
%$$
%that clearly implies
%$$
%\gM(S_\vf)=\sum_{i=1}^N | k_i^+-k_i^- |\,\cH^2(U_i)\,.
%%\quad\fa\, (k_1,\ldots,k_N)\in \text{{\rm Ad}}_{\min}(\vf) \,.
%$$
Since $k_i^\pm\in\gn$, we have $| k_i^+ - k_i^- |
{\leq}
k_i^+ + k_i^- =k_i$ for every $i$ and hence the inequality \eqref{ineq} readily follows.

If $[\vf]$ satisfies property (H), by Corollary \ref{cor}  equation \eqref{nogap} is satisfied.
Assume now that $[\vf]$ does not satisfy property (H), so we have to show that the inequality in \eqref{ineq} is strict. We can  then find an index $i\in\{1,\ldots,K\}$ for which $|  \#(\s_i^+) - \#(\s_i^-) |\neq k_i=k_i^++k_i^-$. Recalling that
$\#(\s_i^+),\,k_i^\pm\in\gn$, we thus have $|  \#(\s_i^+) - \#(\s_i^-) |< k_i$ where $\#(\s_i^+) - \#(\s_i^-) =k_i^--k_i^+$ by \eqref{diff}.
Therefore, $| k_i^+-k_i^- |< k_i$ for some $i$, whereas in general $| k_i^+-k_i^- |\leq k_i$.
In conclusion, the strict inequality $\gM(S_\vf)<P(\vf)$ holds, as required. \epf
%
%%%%%%%%
%
\section{{Relaxation in $L^1$}}\label{sec:5}
{We first turn our attention to the relaxation of $\mathcal F_f$ in $L^1$. We focus on the special case of the vortex map
$$u_V\in W^{1,1}(B^n_R;\gR^n)\,,\quad u_V(x):=\frac{x}{|x|}\,, $$
where we consider as domain an open ball $B^n_R$ in $\Rn$ centered at the origin $0_{\gR^n}$ and of radius $R>0$.

We first discuss the relaxed area in dimension $n=2$.
Following \cite{BES1,BES2,BES3}, we can consider the localized functional $\wid\cA_{L^1}(u_V,B^2_R)$ for $R$ small, for which it happens that the energy gap of the relaxed area is equal to the mass of a suitable vertical current $S_T\in\cR_2(B^2_R\tim\gR^2)$ that is supported on a catenoidal surface $\Sigma_R$ whose projection onto the domain lies in a segment $I_R$ connecting the origin $0_{\gR^2}$ to a given point $P_R$ at the boundary $\pa B^2_R$. Such a surface $\Sigma_R$ is obtained through a minimization problem for the surface area in a $3$-dimensional Plateau type problem in $I_R\tim\gR^2$ with a free boundary in correspondence to $P_R$ (see \cite{BES2}). Instead, if $R>0$ is big, it is easy to see (this was already observed in \cite{AcDM}) that the energy gap is equal to $\pi$, namely
\begin{equation}
\wid\cA_{L^1}(u_V,B^2_R)=\mathcal V(u_V,B^2_R)+\pi=\mathcal V(u_V,B^2_R)+ P \Bigl(\frac{x}{|x|} \Bigr),
\end{equation}
where $P(\frac{x}{|x|})=\pi$ is in fact the mass of the minimal current $S_V$ solution of the disk type Plateau problem spanning the unit circumference $\mathbb S^1$. In this case, a recovery sequence $u_h\in C^1(B^2_R;\gR^2)$ for $\wid\cA_{L^1}(u_V,B^2_R)$ is such that {(as currents in $\cD_2(B^2_R\times\gR^2)$)}
$$
G_{u_h}\wc T_{u_V}= G_{u_V}+\delta_{0_{\gR^2}}\times \qu{D^2}\,,\quad
\gM(G_{u_h})\rightarrow \gM(T_{u_V})\,,
$$
{where $D^2$ is the unit disk centered at the origin in the target space.}

Let now $f$ be a continuous integrand satisfying hypotheses
$(a),\,(b),\,(c)$ from Sec. \ref{Sec:poly}.
Let $(v_h)\subset C^1(B^2_R;\gR^2)$ be a sequence tending to $u_V$ in $L^1(B_R^2;\gR^2)$ and such that  $\cE_f(v_h,B^2_R)\to \wid{\mathcal{F}}_{L^1}(u,B^2_R)$ as $h\to\infty$. Up to a subsequence, we can suppose that $G_{v_h}\rightharpoonup T\in \mathcal T_{u_V}$ weakly as currents. If we show that
\begin{equation}\label{44}
\left\{\ba{l}T=T_{u_V}= G_{u_V}+\delta_{0_{\gR^2}}\times \qu{D^2} \\
\lim\limits_{h\to \infty}\GM(G_{v_h})=\GM(T_{u_V})\,,
\ea
\right.
\end{equation}
we could apply Proposition \ref{Pcont} to  deduce that
\begin{equation}
\wid{\mathcal{F}}_{L^1}(u_V,B^2_R)=\int_{B^2_R}f(x,u_V(x),\nabla u_V(x))\dd x+\mathcal F_f(\delta_{0_{\gR^2}}\times \qu{D^2}).
\end{equation}
{For integrands $f$ not depending on $x$, Proposition \ref{Pvortex}, or under an additional structure hypothesis on $f$, Proposition \ref{Pvortex2}, we} show that \eqref{44} holds if $R$ is large enough.
\subsection{Estimates on the energy gap}
For $n=N\geq 2$, we thus consider a non-negative continuous integrand $f$
satisfying hypotheses
$(a),\,(b),\,(c)$ from Sec. \ref{Sec:poly}, so that we have
$$ c\,|M(G)|\leq f(x,u,G)\leq C\,|M(G)|\,, $$
and we consider the localized relaxed energy $\wid\cF_{L^1}(u_V,B^n_R)$ of the vortex map $u_V:B^n_R\to\Rn$.

In the sequel, we let $C^n_f$ denote the $\cF_f$-energy of the vertical current
$\d_{0}\tim\qu{D^n}$, where $0=0_{\gR^n}$ and $D^n$ is the unit $n$-disk centered at the origin in the target space, that is
$$ C^n_f=\cF_f( \d_{0}\tim\qu{D^n})\,.
$$
Notice that we have $c\,\o_n\leq C^n_f\leq C\,\o_n$, where we denote $\o_n=|D^n|$.
\bp\label{Pvortex} Let $f$ be a continuous integrand satisfying properties (a), (b), (c), and independent of $x$, $f(x,u,G)=\widehat f(u,G)$. There exists a positive radius $R_0=R_0(n,c,C)$ such that if $R>R_0$ we have
\begin{equation}
\label{vortex-en}
\wid\cF_{L^1}(u_V,B^n_R) = \cE_f(u_V,B^n_R)+C^n_f\,.
\end{equation}
\ep
\bpf For every $R>0$, the relaxed energy on $B^n_R$ is lower than the energy of the Cartesian current $G_{u_V}+\d_0\tim\qu{D^n}$ on the cylinder $W_R:=B^n_R\tim\Rn$, that is equal to $\cE_f(u_V,B^n_R)+C^n_f$ (this is true since there exists a sequence $v_h$ for which \eqref{44} holds). This yields the inequality "$\leq$".

To prove the opposite inequality, for any $R>0$ we let $\{u_h\}\sb C^1(B^n_R;\RN)$ be such that $u_h\to u_V$ in $L^1(B^n_R;\Rn)$ and $\cE_f(u_h,B^n_R)\to \wid\cF_{L^1}(u_V,B^n_R)$. We can also assume that $G_{u_h}$ weakly converges to some current $T=G_{u_V}+T_s$ in $\cart(B^n_R\tim\Rn)$. By lower semicontinuity, it suffices to show that if $R$ is large enough we have
\begin{equation}
\label{Cflb}
\cF_f(T_s)\geq C^n_f\,.
\end{equation}
To this purpose, arguing as in \cite{AcDM} we obtain the following
\bl\label{lem:slicing} If $R$ is large enough, the projection on $B^n_R$ of the support of $T_s$ is a compact subset of $B^n_R$.
\el
\bpf Following exactly the same slicing argument as in \cite[Lemma 5.2]{AcDM}, we find that the function $\z(t)$ satisfies
$$ 0\leq \lim_{t\to R^-}[\z(t)]^{1/n}\leq[\z(0)]^{1/n}-\frac R{n\,\g_n^{1-1/n}}, $$
where $\gamma_n$ is the isoperimetric constant.
On the other hand, using that $\cF_f(T_s)\leq C^n_f$, we can estimate
$$
\z(0)\leq \GM_{W_R}(T_s)\leq \frac 1c\,\cF_f(T_s)\leq  \frac 1c\,C^n_f\leq \frac C c\,\o_n\,,
$$
and hence we obtain a contradiction, if $R$ is sufficiently large.
\epf

Now, the previous lemma implies that if $R$ is large enough, we can see $T_s$ as an integral current in $\Rn\tim\Rn$ satisfying $\pa T_s =\d_{0}\tim\qu{\Sph^{n-1}}$, where $\Sph^{n-1}=\pa D^n$. {We claim that the energy of the current $T_s$ is greater than the energy of its image} through the orthogonal projection $\Pi(x,y)=(0,y)$.
Since
$$ \pa(\Pi_\sharp T_s)=\Pi_\sharp(\pa T_s)=\Pi_\sharp(\d_0\tim\qu{\Sph^{n-1}})=\d_0\tim\qu{\Sph^{n-1}}\,, $$
by the Constancy theorem we have
$$\Pi_\sharp T_s= \d_{0}\tim\qu{D^n}  $$
and hence we estimate
$$
\cF_f(T_s)\geq \cF_f(\Pi_\sharp T_s) = \cF_f(\d_0\tim\qu{D^2}))\,,
$$
so that \eqref{Cflb} holds.
{To prove the claim, we recall that since the integrand $f$ does not depend on $x$, the integral formula \eqref{7pF} of the $\cF_f$-energy gives
$$ \cF_f(T_s)=\int_{B_R^n\times \gR^n} F_f\bigl(\wih\p(z),\vect
{T_s}(z)\bigr)\,\dd\Vert T_s\Vert(z)\,,\quad
\cF_f(\Pi_\sharp T_s)=\int_{B_R^n\times \gR^n} F_f\bigl(\wih\p(z),\Lambda_n\Pi(\vect
{T_s}(z))\bigr)\,\dd\Vert \Pi_\sharp T_s\Vert(z)\,, $$
where clearly $\Vert T_s\Vert\geq\Vert\Pi_\sharp T_s\Vert$. Furthermore, for every $u=\wih\p(z)$ the integrand
$F_f\bigl(u,\cdot\bigr)$ defines a seminorm on the vector space $\Lambda_n\gR^{n+n}$, so that for every $\xi\in \Lambda_n\gR^{n+n}$ we have $F_f\bigl(u,\xi\bigr)\geq F_f\bigl(u,\Lambda_n\Pi(\xi)\bigr)$. Therefore, the claim follows and the proof is complete.
}
\epf

\par We now observe that the {localized} functional $A\mapsto\wid{\mathcal{F}}_{L^1}(u_V,A)$ fails to be subadditive, and hence it cannot be extended to a measure on Borel {subsets of $\Rn$.}

This non-locality property is checked by means of the very same argument as in \cite{AcDM}, on account of Proposition \ref{Pvortex} and of the following upper estimate.
\bp
\label{Pvortex2}
For every $R>0$ we have
$$ \wid\cF_{L^1}(u_V,B^n_R) \leq \cE_f(u_V,B^n_R)+C\,\frac{\o_n}n\,R\,.
$$
\ep
\bpf {Let $T_R$ be the current
$$ T_R=G_{u_V}+L_R\tim\qu{\Sph^{n-1}} \in\cart(B^n_R\tim\Rn),
$$
where $L_R$ is the $1$-current integration on a line segment connecting a point at the boundary of $B^n_R$ with the origin.
We can easily find a sequence $\{u_h\}\sb C^1(B^n_R;\Rn)$ converging to $u_V$ in $L^1$ and such that $G_{u_h}\wc T_R$. Therefore, the relaxed energy of $u_V$ is lower than the energy of the current $T_R$.} 
We have
$$ \cF_f(T_R)=\cE_f(u,B^n_R)+\cF_f(L_R\tim\qu{\Sph^{n-1}})
$$
where by the upper estimate
$$  \cF_f(L_R\tim\qu{\Sph^{n-1}})\leq C\, \GM(L_R\tim\qu{\Sph^{n-1}})=C\,R\,\cH^{n-1}(\Sph^{n-1})\,,
$$
and the assertion readily follows, again arguing as in \cite[Lemma 5.2]{AcDM}.
\epf

In the case of integrands $f$ which also depend on the {$x$ variable} we need an additional {hypothesis} to conclude as in Proposition \ref{Pvortex}.
{Below, we denote for simplicity $\vect\e:=\e_1\wedge\cdots\wedge\e_n$, the unit simple $n$-vector orienting the target space $\Rn$ as a subspace of $\gR^{n+n}$.}

\bp
Assume that 
\begin{equation}\label{46}
	{F_f(x,u,\xi)\geq \max_{u\in D^n}\{F_f(0_{\gR^n},u,\vect\varepsilon)\}=:\widetilde c,}\qquad \qquad \forall\, x,u\in \gR^n,\;\xi\in \Lambda_n\gR^{n+n}\; \text{with }|\xi|=1.
\end{equation}
Then there exists a constant $\overline R>0$ such that for all $R\geq \overline R$ {property \eqref{vortex-en} holds.}
%\begin{equation}
%	\wid{\mathcal{F}}_{L^1}(u_V,B^n_R)=\mathcal E_f(u_V,B^n_R)+C^n_f.
%\end{equation}
\ep

\bpf
Arguing as above, Lemma \ref{lem:slicing} implies that for $R\geq \overline R$ there holds
%\begin{equation}
$$
\mathcal F_f(T_s)=\int_{\Sigma} |\theta(x,u)|  F_f(x,u,\tau(x,u))\,\dd\mathcal H^{n}(x,u)\geq \widetilde c\int_{\Sigma} |\theta| \,\dd\mathcal H^{n}=\widetilde c\,\gM(T_s)\geq \widetilde c\,\omega_n,
$$
%\end{equation}
where we have written $T_s=\qu{\Sigma,\theta,\tau}$ and used the fact that, since $\partial T_s=\delta_{0}\times\qu{\mathbb S^{n-1}}$, then $\gM(T_s)\geq \gM(\delta_{0}\times \qu{D^n})=\omega_n$.
Therefore, again by \eqref{46}
%\begin{equation}\label{47}
$$
\mathcal F_f(T_s)\geq \widetilde c\,\omega_n\geq \int_{\{0\}\times D^n}F_f(0,u,\vect\varepsilon)\,\dd\mathcal H^{n}=C^n_f,
$$
%\end{equation}
hence $\wid\cF_{L^1}(u,B^n_R)\geq\mathcal E_f(u_V,B^n_R)+C^n_f.$

Let us now show the opposite inequality: since for $\wid\cA_{L^1}(u,B^n_R)$ there is always a sequence $\{u_k\}\subset C^1(B^n_R;\Rn)$ such that $G_{u_k}\wc G_{u_V}+\delta_{0}\times \qu{D^n}$ and  $A(u_k,B^n_R)\rightarrow  \gM(G_{u_V})+\gM(\delta_{0_{\gR^2}}\times \qu{D^2})$, we readily infer by Proposition \ref{Pcont}  that $\wid\cF_{L^1}(u,B_R)\leq\mathcal E_f(u_V,B^n_R)+C^n_f.$
\epf
%
%
%
%, it is easy to see that if $R>0$ is big enough, then $\delta_{0_{\gR^2}}\times \qu{D^2}$ is the minimal (in mass) integral current $\widehat T$ in $\mathcal R_2(B_R\times \gR^2)$ such that $\partial \widehat T=-\partial G_{u_V}$. If we knewn that this is minimal also for $\mathcal F_f$ we could apply Proposition \ref{Pcont} to  . Precisely,
%assume also that
%\begin{equation}
%\mathcal F_f(\widehat T)\leq \mathcal F_f(T_v)\;\qquad \forall T_v\in \mathcal R_2(B_R\times \gR^2)\;\text{ with }\;\partial T_v=-\partial G_{u_V},
%\end{equation}
%then Corollary \ref{cor4.7} implies that
%$$\wid\cF_{L^1}(u,B^n)=\cF_f( G_{u_V}+\widehat T)=\cE_f(u_V,B_R)+\cF_f(\delta_{0_{\gR^2}}\tim\qu{D^2})\,. $$}
%
\section{Some explicit formulas}\label{sec:6}
We now discuss some situations where Theorem \ref{T-cont} applies: starting from known results for the relaxed area functional {in the strict convergence}, we obtain explicit formulas for more general energies as above.
\subsection{Sobolev maps into the unit circle}\label{sec:sobolev}
Denote by
$$\Sph^1:=\{y\in\gR^2\,:\,|y|=1\}\,,\quad D^2:=\{y\in\gR^2\,:\,|y|<1\} $$
the unit circle and disk in the target space, and  let
$$ W^{1,1}(B^n;\Sph^1):=\{u\in W^{1,1}(B^n;\gR^2): |u(x)|=1\quad{\mx{\text for $\calL^n$-a.e. }}x\in B^n\}\,. $$
\par In any dimension $n\geq 2$, the graph current of a Sobolev map
$u\in W^{1,1}(B^n;\Sph^1)$ is an i.m. rectifiable $n$-current $G_u$ in $B^n\tim\Sph^1$, with finite mass
$$ \gM(G_u)=\int_{B^n}\sqrt{1+|\nabla u|^2}\,\dd x=\cV(u,B^n)<\i\,. $$

In low dimension $n=2$, the following result was obtained in \cite{BSS}. It says that the energy gap is detected by the {\em distributional determinant} $\Det\nabla u$.
\bt\label{T-BSS} Let $u\in W^{1,1}(B^2;\Sph^1)$, and let $\Det\nabla u$ denote the {\em distributional determinant} of $u$. Then,
$$ \wid{\mathcal{A}}_{BV}(u,B^2)<\i \iff |\Det\nabla u|(B^2)<\i\,. $$
In that case, moreover, one has:
$$ \wid{\mathcal{A}}_{BV}(u,B^2) =\cV(u,B^2) +|\Det\nabla u|(B^2)\,. $$
\et
\par In any dimension $n\geq 2$, the relevant singularities of $u$ are detected by the current $G_u$, i.e., they can be described by homological arguments.
More precisely, denoting by {$\o_{\Sph^1}$} the closed 1-form in $\Sph^1$
$$ %\beq\label{o2}
{\o_{\Sph^1}}:=\ds \frac 12\,\bigl(y^1\,dy^2-y^2\,dy^1 \bigr) $$
the singularities are read by the $(n-2)$-dimensional current $\rP(u)\in\cD_{n-2}(B^n)$ defined by
$$ \rP(u)(\y):= -\frac 1\pi G_u(\dd\y\wedge{\o_{\Sph^1}})\,,\quad \y\in \cD^{n-2}(B^n)\,. $$
Therefore, in low dimension $n=2$ one has
\beq\label{DetPu}
\pi\cdot\rP(u)(\y)=\langle\Det\nabla u,\y\rangle\quad\fa\,\y\in C^\i_c(B^2)\,.
\eeq
If e.g. $u_V(x)=x/|x|$, the vortex map, one gets $\rP(u_V)=\d_{0}$, the unit Dirac mass at the origin.

Theorem \ref{T-BSS} was extended in \cite{CM} as follows:
\bt\label{T-CM}
Let $n\geq 2$ and $u\in W^{1,1}(B^n;\Sph^1)$. Then, $\wid{\mathcal{A}}_{BV}(u,B^n)<\i$ if and only if the $(n-2)$-current $\rP(u)$ is i.m. rectifiable and with finite mass, $\gM(\rP(u))<\i$.
In that case, moreover, one has:
$$ \wid{\mathcal{A}}_{BV}(u,B^n) =\cV(u,B^n)+\pi\,\gM(\rP(u))\,. $$
\et
\par We recall the strategy of the proof. A result from \cite{Mu22} implies that if $u\in W^{1,1}(B^n;\Sph^1)$ satisfies $\wid{\mathcal{A}}_{BV}(u,B^n)<\i$, then there exists a unique optimal Cartesian current $T_u\in \mathcal T_u^{BV}$ that encloses the graph of $u$, and it is given by
$$ T_u=G_u+(-1)^{n-2}\rP(u)\times\qu{D^2}\,. $$
Therefore, the proof of the energy lower bound readily follows: by Federer-Fleming Compactness theorem, for every smooth sequence $\{u_h\}$ strictly converging to $u$, the graph $G_{u_h}$ weakly converges to $T_u$ (up to extracting a subsequence), and one concludes from the semicontinuity of the mass.
On the other hand, the energy upper bound holds true as a consequence of the following approximation result from \cite{CM}:
\bt\label{T-CMappr}
Let $n\geq 2$ and $u\in W^{1,1}(B^n;\Sph^1)$ be a Sobolev map with finite relaxed energy, $\wid{\mathcal{A}}_{BV}(u,B^n)<\i$. Then, there exists a sequence
$\{u_h\}\sb C^\infty(B^n;\gR^2)$ such that $G_{u_h}\wc T_u$ weakly in $\cD_n(B^n\times\gR^2)$ and
$\gM(G_{u_h})\to \gM(T_u)$ as $h\to\infty$.
\et
\par As a consequence of Theorems \ref{T-cont} and \ref{T-CMappr}, we readily obtain the following density result.
\bp\label{P-CMappr}
Let $n\geq 2$ and $u\in W^{1,1}(B^n;\Sph^1)$ be a Sobolev map with finite relaxed energy, $\wid{\mathcal{A}}_{BV}{(u,B^n)}<\i$. Let $f$ be a continuous integrand satisfying properties
{$(a)$, $(b)$ from Sec. \ref{Sec:poly}, where $N=2$.} Then, there exists a smooth sequence
$\{u_h\}\sb C^\infty(B^n;\gR^2)$ such that $G_{u_h}\wc T_u$ weakly in $\cD_n(B^n\times\gR^2)$ and
$$\lim_{h\to \i}\cE_f(u_h,B^n)=\cE_f(u,B^n)+\cF_f(\rP(u)\tim\qu{D^2})\,. $$
\ep

\par Arguing as above, and using the lower semicontinuity result in Proposition \ref{Plsc},
we thus readily obtain:
\bt\label{T-S1}
Let $n\geq 2$ and let $f$ be a continuous integrand satisfying properties
{$(a),\,(b),\,(c')$ in Sec. \ref{Sec:poly},}
where $N=2$. For any $u\in W^{1,1}(B^n;\Sph^1)$, we have $\wid{\mathcal{F}}_{BV}(u,B^n)<\i$ if and only if the $(n-2)$-current $\rP(u)$ is i.m. rectifiable and with finite mass, $\gM(\rP(u))<\i$. In that case, moreover, one has:
$$ \wid{\mathcal{F}}_{BV}(u,B^n) =\cE_f(u,B^n)+\cF_f(\rP(u)\tim\qu{D^2})\,. $$
\et

We would like to find an explicit formula of the {\em energy gap} obtained in Theorem \ref{T-S1}, i.e., of the integral
\beq\label{gap} \cF_f(\rP(u)\tim\qu{D^2})\,.
\eeq
We first consider the easier case when $f=f(G)$
{is a polyconvex function
%does not depend on $(x,u)$, whence $f:\gR^{N\tim n}\to\gR_+$ is convex and lower
%semicontinuous, and
satisfying properties $(b)$ and $(c')$, that is
$$ c\,|M_2(G)|\leq f(G)\leq C\,|M(G)| \quad \fa\, G\in \gR^{2\tim n}. $$
%where we let $N=2$.
%
Therefore, for smooth functions $u:B^n\to\gR^2$ we have}
%, eq. \eqref{Ef} becomes
$$ \cE_f(u,B^n):=\int_{B^n}f(\nabla u(x))\,\dd x\,.
$$
More generally, if $T\in\cR_n(B^n\tim\gR^2)$, and $T=\qu{S,\t,\x}$, we have
$$ \cF_f(T)= \int_{S} {\t(z)}\,F_f\bigl(\x(z)\bigr)\,d\cH^n(z)\,. $$

Notice that if $u\in W^{1,1}(B^n;\Sph^1)$ is such that  $\wid{\mathcal{F}}_{BV}{(u,B^n)}<\i$, in high dimension $n\geq 3$ the energy gap in \eqref{gap} depends on the orienting $(n-2)$-vector to the current $\rP(u)$. However, in low dimension $n=2$, recalling that the 2-current $\qu{D^2}$ is oriented by $\e_1\wedge\e_2$, and using that $\rP(u)$ is 0-dimensional, by eq. \eqref{DetPu} we obtain
$$ \cF_f(\rP(u)\tim\qu{D^2})=\GM(\rP(u))\cdot \int_{D^2} F_f\bigl(\e_1\wedge\e_2\bigr)\,d\cH^2(y) = |\Det\nabla u|(B^2)\cdot F_f\bigl(\e_1\wedge\e_2\bigr)\,.
$$
In fact, if
{$f(G)=|M(G)|$ or $f(G)=|\det G|$}, we clearly have $F_f\bigl(\e_1\wedge\e_2\bigr)=1$, see Theorem \ref{T-BSS}.
%
%%%%%%
%
\subsection{The case of $0$-homogeneous maps}
Similar arguments apply to the case of homogeneous $BV$-maps treated in \cite{Ca24}. Let $\gamma\in BV(\mathbb S^1;\gR^2)$ be fixed, and consider the function $u_\gamma\in BV(B^2;\gR^2)$ given by
\begin{equation}\label{ug}
u_\gamma(x)=\gamma {\left(\frac{x}{|x|} \right)},\qquad \qquad x\in B^2\setminus \{0\}.
\end{equation}
In \cite{Ca24} it is proved that
$$\wid{\mathcal{A}}_{BV}(u_\gamma,B^2)=\mathcal V(u_\gamma,B^2)+\overline P(\gamma),$$
where, according to definition in \eqref{plateau_rel} and \eqref{plateau_rel=plateau} one has
$$\overline P(\gamma)=P(\overline \gamma)$$
where $\overline \gamma$ is the Lipschitz curve in Definition \ref{defrep}. If we assume that $\overline \gamma$ satisfy hypothesis (P0), then we obtain the following result:

\begin{theorem}\label{teo_0hom1}
Let $\gamma\in \text{{\rm Lip}}(\mathbb S^1;\gR^2)$  satisfy hypothesis $(P0)$, let $u_\gamma\in W^{1,1}(B^2;\gR^2)$ be defined as in \eqref{ug}, and let $f$ be a continuous integrand satisfying properties
{$(a),\,(b),\,(c')$ from Sec. \ref{Sec:poly},}
where $N=2$.
If $\gamma$ {satisfies property $(H)$} in Definition \ref{defH} then
\begin{equation}
\wid{\mathcal{F}}_{BV}({u_\gamma},B^2)=\cE_f({u_\gamma},B^2)+	P_f(0,\gamma),
\end{equation}
where
\begin{equation}
\label{Pfg}
P_f(0,\gamma)=\sum_{i=1}^{K}k_i\int_{U_i}F_f(0_{\gR^2},y,\alpha_i\varepsilon_1\wedge\varepsilon_2)dy
\end{equation}
with $(k_1,\dots,k_K)$ is the unique injection in $\text{{\rm Ad}}_{\min}(\gamma),$ and where $S_\gamma= \sum_{i=1}^K \alpha_i k_i\,\qu{U_i}$, $U_i$ being the
{bounded} connected components of $\gR^2\setminus \gamma(\mathbb S^1)$ and $\alpha_i$ suitable signs (see Theorem \ref{teo2.9}).
\end{theorem}

\bpf
The completely vertical component $S_T\in \mathcal R_2(B^2\times \gR^2)$ of the unique current $T_{u_\gamma}=G_{u_\gamma}+S_T\in \mathcal T_{u_\gamma}^{BV}$, is given by $$S_T=\delta_{0}\times   S_{\gamma},$$ and since $\gamma$ satisfies (H), by Theorem \ref{Tgap} it holds
$$\gM(S_T)=\gM(S_\gamma)=P(\gamma).$$
Hence the thesis follows from Theorem \ref{T-cont}, by observing that $S_\gamma$, by Theorem \ref{teo2.9}, has the form
$$S_\gamma= \sum_{i=1}^K \alpha_i k_i\,\qu{U_i}, \qquad \qquad \alpha_i\in \{-1,1\},$$
and $\alpha_i\,\varepsilon_1\wedge\varepsilon_2$ orients $S_T$ on $\{0\}\times U_i$.
\epf

We have already seen that the previous result does not apply in the case of the Lipschitz function {$\vf_8$ in \eqref{ottodoppio}} corresponding to the double eight example {$u_8$}, for which property (H) fails to hold.
However, we have
\bp\label{prop:8}  Let $f$ be a continuous integrand satisfying properties
$(a),\,(b),\,(c)$ from Sec. \ref{Sec:poly}, and let ${u_8}$ be the homogeneous map given by the double eight example. Then, the {energy} gap is positive.
\ep
\bpf Assume that we have
$$ \wid{\mathcal{F}}_{BV}({u_8},B^2)=\cE_f({u_8},B^2)\,. $$
Then, by the definition of relaxed functional it turns out that for every small radius $R>0$ we have
$$ \wid{\mathcal{F}}_{BV}({u_8},B^2_R)=\cE_f({u_8},B^2_R)\,. $$
 {By the estimate $(c)$ we can find a constant $c>0$ such that}
$$ \wid{\mathcal{A}}_{BV}({u_8},B^2_R)\leq \frac 1c\,\wid{\mathcal{F}}_{BV}({u_8},B^2_R)
=\frac 1c\,\cE_f({u_8},B^2_R)\,,$$
where the right-hand side is small if $R$ is small, by absolute continuity. However, it is well-known that in the double eight example we have
$$ \lim_{R\to 0^+}\wid{\mathcal{A}}_{BV}({u_8},B^2_R)=2\pi\,, $$
and hence we obtain a contradiction.
\epf
\subsubsection{The BV case}
We now observe that, if $\gamma$ is not Lipschitz but only of bounded variation, {in} the previous discussion $P(\gamma)$ must be replaced by $P(\overline \gamma)$. In such a case, however, some additional vertical part in $S_T$ pops up: specifically, let us
{view $\g$ as a function of the angle variable $\theta$ and} denote by $J_\gamma$ the jump set of $\gamma$, $J_\gamma=\{\theta_\ell: \ell\in \mathbb N\}\subset \mathbb S^1$.
{Without loss of generality, we assume here that the jump set of $\g$ is countable and that $\g$ is continuous at $(0,1)$.
It turns out that}
\begin{equation}
D_\theta\gamma=D_\theta^a\gamma+D^C_\theta \gamma+D_\theta^J\gamma=\dot\gamma \cdot\mathcal H^1+D^C_\theta \gamma+\sum_{\ell}(\gamma^+(\theta_\ell)-\gamma^-(\theta_\ell))\,\delta_{\theta_\ell}.
\end{equation}
The three measures above are mutually singular, and we denote by $C'_\gamma\subset\mathbb S^1$ a measurable set of null {$\Ha^1$} measure such that Cantor part $D^C\gamma$ satisfies $|D^C\gamma|(\mathbb S^1\setminus C'_\gamma)=0$. By identifying $\mathbb S^1$ with $[0,2\pi)$ and using polar coordinate in $B^2$, we set
$C_\gamma:=\{x\in B^2\setminus \{0\}:x=(\rho\cos\theta,\rho\sin\theta),\;\theta\in C'_\gamma\}$; we
can easily see that $Du_\gamma=\nabla u_\gamma+D^C u_\gamma+D^Ju_\gamma$, with
the unique current $T_{u_\g}\in\cT^{BV}_{u_\g}$ given by
\begin{equation}
T_{u_\gamma}=G_{u_\g}+S_T,\qquad \qquad S_T=S_\g^C+S_\g^J+\delta_{0_{\gR^2}}\times   S_{\ol\gamma}\,,
\end{equation}
where $\ol\g$ is the Lipschitz curve from Definition \ref{defrep}.

The components $G_{u_\g}$, $S_\g^C$ and $S_\g^J$ are respectively given by
$$ G_{u_\g}=\Psi_\#(\qu{0,1}\tim T^a_\g)\,,\quad S_\g^C=\Psi_\#(\qu{0,1}\tim T^C_\g)\,,\quad S_\g^J=\Psi_\#(\qu{0,1}\tim T^J_\g)\,, $$
where $\Psi:(0,1)\tim [0,2\pi)\tim\gR^2\to B^2\tim\gR^2$ is the map
$$ \Psi(\r,\t,y)=(\r\cos\t,\r\sin\t,y) $$
and the terms $T^a_\g$, $T^C_\g$, and $T^J_\g$ in $\cR_1((0,2\pi)\tim\gR^2)$ correspond to the absolutely continuous, Cantor and Jump components.
More precisely, $T^a_\g$ is the graph current of the $BV$ function $\g:(0,2\pi)\to\gR^2$. Moreover, the components $T^C_\g$ and $T^J_\g$ are vertical, i.e., for every test function $\f\in C^\i_c(\gR\tim\gR^2)$ we have
$$ {T^C_\g( \f(\t,y)\,d\t )} = { T^J_\g( \f(\t,y)\,d\t )}=0,$$
whereas for $j=1,2$, denoting by components $\g=(\g_1,\g_2)$, we have
$${ T^C_\g ( \f(\t,y)\,dy^j) } =\int_{C'_\g}\f(\t,\g(\t))\,dD^C {\g_j}   $$
and
$$ { T^J_\g (\f(\t,y)\,dy^j) } =\sum_\ell\int_0^1\f(\t_\ell,s\,\g^+(\t_\ell)+(1-s)\,\g^-(\t_\ell))\,( \g_j^+(\t_\ell)-\g_j^-(\t_\ell) )\,ds \,.$$

Then, for every test function $\phi\in C^\infty_c(B^2\times \gR^2)$, we have
$$ { S_\g^C(\f(x,y)\,dx^1\wedge dx^2)}={ S_\g^J(\f(x,y)\,dx^1\wedge dx^2)}=0 $$
and for $i,j\in\{1,2\}$
$$ { S_\g^C(\f(x,y)\,\wih{\dd x^i}\wedge\dd y^j )} =(-1)^{i-1}\int_{B^2} \f(x,u(x))\,d(D^C u_\g)^j_i \,,$$
where $u_\g$ is a precise representative, and finally
\begin{equation}\label{58}
{ S_\g^J (\f(x,y)\,\wih{\dd x^i}\wedge\dd y^j) }=(-1)^{i-1}\sum_\ell\int_{R_\ell}\left(\int_0^1\phi(x,u^s(x))\dd s\right)(\gamma^+_j(\theta_\ell)-\gamma^-_j(\theta_\ell)){\nu^i_\ell}\,\dd \mathcal \cH^1 \,,
\end{equation}
where $R_\ell$ is the outward pointing radius of $B^2$ oriented by $(\cos\theta_\ell,\sin\theta_\ell)$, the unit normal is $\nu_\ell=(\nu_\ell^1,\nu_\ell^2)=(-\sin\theta_\ell,\cos\theta_\ell)$, and for any $x\in R_\ell$ and any $\ell\in\gn$ we have denoted
$$ u^s(x)=s\,\g^+(\t_\ell)+(1-s)\,\g^-(\t_\ell)\,. $$

We thus have
$$ \GM(T_{u_\g})=\GM(G_{u_\g})+\GM(S^C_\g)+\GM(S^J_\g)+\GM(S_{\ol\g})\,, $$
where
$$ \GM(G_{u_\g})=A(u_\g,B^2)\,,\quad \GM(S^C_\g)=|D^C u_\g|(B^2)\,,\quad
\GM(S^J_\g)=\sum_\ell |\gamma^+(\theta_\ell)-\gamma^-(\theta_\ell) |\,.
$$
Therefore, using the results from \cite{Ca24} it turns out that
$$\wid\cA_{BV}(u_\gamma,B^2)=\GM(T_{u_\g})\,. $$

By applying Theorem \ref{T-cont}, we thus readily obtain the following representation result:
\begin{theorem}\label{teo_0hom2}
	Let $\gamma\in BV(\mathbb S^1;\gR^2)$ be such that $\overline \gamma$ satisfy hypothesis (P0), let $u_\gamma\in BV(B^2;\gR^2)$ be defined as in \eqref{ug}, and let $f$ be a continuous integrand satisfying properties
{$(a),\,(b),\,(c)$ from Sec. \ref{Sec:poly},}
%$(a),\,(b),\,(c)$,
where $N=2$.
	If $\overline \gamma$ satisfies hypothesis $(H)$ in Definition \ref{defH} then
	\begin{equation}
		\wid{\mathcal{F}}_{BV}(u_\g,B^2)=\cE_f(u_\g,B^2)+\cF_f(S^C_\g)+\cF_f(S^J_\g)+P_f(0,\ol\gamma)
\end{equation}
where $P_f(0,\ol\gamma)$ is defined as in \eqref{Pfg} with respect to the Lipschitz curve $\ol\g$ from Definition \ref{defrep}.
	\end{theorem}
	
	\textbf{An example: the triple junction map:}
Consider the function $\gamma:[0,2\pi)\rightarrow \gR^2$ defined as
$$\gamma(\theta):=\begin{cases}
\alpha&\text{ if }\theta\in [0,\pi/3)\cup[5\pi/3,2\pi)\\
\beta&\text{ if }\theta\in [\pi/3,\pi)\\
\gamma&\text{ if }\theta \in [\pi,5\pi/3),
\end{cases}$$
where $\alpha,\beta,\gamma\in \mathbb S^1$ are the vertices of an equilateral triangle $\Delta$; e.g., $\alpha=(1,0)$, $\beta=(-\frac12,\frac{\sqrt3}{2})$, $\gamma=(-\frac12,-\frac{\sqrt3}{2})$.
Consider the triple junction map  $u_T:B^2\rightarrow \gR^2$ given by
\begin{equation}\label{triplemap}
	u_T(\rho,\theta)=\gamma(\theta).
\end{equation}
 Denoting by $\overline{pq}$ the oriented segment from $p\in \mathbb R^2$ to $q\in \gR^2$, we have $T_\gamma^a=G_{\gamma}$, $T_\gamma^C=0$, whereas $T_\gamma^J=\delta_{\pi/3}\times \qu{\overline{\alpha\beta}}+\delta_{\pi}\times \qu{\overline{\beta\gamma}}+\delta_{5\pi/3}\times \qu{\overline{\gamma\alpha}}$, and so $S_\gamma^C=0$.
Moreover \eqref{58} holds with $\ell\in \{1,2,3\}$ and with $\theta_1=\pi/3$, $\theta_2=\pi$, $\theta_3=5\pi/3$, $\nu_1=(-\sqrt3/2,1/2)$, $\nu_2=(0,-1)$, {$\nu_3=(\sqrt3/2,1/2)$}, $R_1$, $R_2$ and $R_3$ are the segment on which $u_T$ jumps, i.e. the radii of $B^2$ corresponding to $\theta_1=\pi/3$, $\theta_2=\pi$, $\theta_3=5\pi/3$. Finally $S_{\overline \gamma}=\delta_{0_{\gR^2}}\times \qu{\Delta}$, and it is easy to see that hypotheses of Theorem \ref{teo_0hom2} hold, so that we infer
\begin{equation}
	\wid{\mathcal{F}}_{BV}(u_T,B^2)=\cE_f(u_T,B^2)+\cF_f(S^J_\g)+P_f(0,\ol\gamma)
\end{equation}
where
$$P_f(0,\ol\gamma)=\int_{\Delta}F_f(0_{\gR^2},y,\varepsilon_1\wedge\varepsilon_2)\,\dd y $$
	and
	$$ \cF_f(S^J_\g)=\sum_{\ell=1}^3\int_{R_\ell}\int_0^1 F_f\left(x,s\,\g^+(\t_\ell)+(1-s)\,\g^-(\t_\ell),-\nu_\ell^\perp\wedge \frac{\g^+(\t_\ell)-\g^-(\t_\ell)}{|\g^+(\t_\ell)-\g^-(\t_\ell)|}\right)\,\dd s\,\dd\mathcal H^1,$$
{with $\nu_\ell^\perp$ the unit vector obtained through a counter-clockwise $\pi/2$-rotation of $\nu_\ell$.}

	We will return on this example in the open problem Section \ref{sec:open}, where we discuss the case of $L^1$-relaxation.

\subsection{Piecewise Lipschitz maps}
We consider in this section piecewise Lipschitz maps, a special class of $BV$ maps $u:B^2\to\gR^2$ analyzed in \cite{BSS24}. They are defined as follows: to start with, we define a Lipschitz partition of $B^2$ as a collection $\{\Om_j,j=1,\dots,N,\;N\in \mathbb N\}$ of mutually disjoint non-empty open subsets of $B^2$ having Lipschitz boundaries, and such that $\cup_{j=1}^N\overline \Om_j=\overline {B^2}$. In such a case we denote $\Sigma:=\cup_{j=1}^N\partial \Om_j\cap B^2$ the interface set of the partition. We will consider Lipschitz {partitions} whose interface is a network, according to the following definition:

\begin{definition}\label{def:net}
	A partition interface $\Sigma$ of a Lipschitz partition $\{\Om_j,j=1,\dots,N\}$ of $B^2$ is said to be a network if the following conditions are satisfied: there exists $n\in \mathbb N$ Lipschitz curves $\alpha_\ell:[a_\ell,b_\ell]\rightarrow \overline{B^2}$, $\ell=1,\dots,n$, such that $$\Sigma=\bigcup_{i=1}^n \overline \alpha_\ell([a_\ell,b_\ell])\cap B^2,$$
	and
	\begin{itemize}
		\item[(n1)] {every $\alpha_\ell$} is injective on $(a_\ell,b_\ell)$ and of class $C^2((a_\ell,b_\ell))$, with $\dot\alpha_\ell\equiv 1$, $\dot\alpha_\ell,\ddot\alpha_\ell\in C^0([a_\ell,b_\ell])$, and $\alpha_\ell((a_\ell,b_\ell))\subset B^2$;
		\item[(n2)]  {$\ell_1\neq \ell_2$} implies $\alpha_{\ell_1}((a_{\ell_1},b_{\ell_1}))\cap\alpha_{\ell_2}((a_{\ell_2},b_{\ell_2}))=\varnothing$;
		\item[(n3)] there exists a finite number of points $\{p_1,\dots,p_m\}\subset {B^2}$, $m\geq 0$, called set of junction points, such that for all $\ell$ with $\alpha_\ell(a_\ell)\neq \alpha_\ell(b_\ell)$ one has $\alpha_\ell(\{a_\ell,b_\ell\})\subseteq \{p_1,\dots,p_m\}\cup\partial B^2$;
		\item[(n4)] if $x\in \alpha_\ell(\{a_\ell,b_\ell\})\cap \partial B^2$, then the curve $\alpha_\ell$ is {transversal} to $\partial B^2$ at $x$;
		\item[(n5)] $\ell_1\neq\ell_2$ implies $\alpha_{\ell_1}(\{a_{\ell_1},b_{\ell_1}\})\cap \alpha_{\ell_2}(\{a_{\ell_2},b_{\ell_2}\})\subseteq \{p_1,\dots,p_m\}$;
		\item[(n6)] for all $i=1,\dots,m$ there exist $r>0$ and $N_i\in \mathbb N$ such that $3\leq N_i\leq N$ and for all $r'\in(0,r)$ the ball $B_{r'}(p_i)\subset B^2$ has nonempty intersection with exactly $N_i$ components $\Om_j$ of the partition $\{\Om_j,j=1,\dots,N\}$.
	\end{itemize}
\end{definition}

Notice that a partition $\{\Om_j,j=1,\dots,N\}$ whose interface is a network {enjoys} the property that every curve $\alpha_j$ might reach the boundary $\partial B^2$, but in such a case, if for example $\alpha_{\ell_1}(a_{\ell_1})\in \partial B^2$, then $\alpha_{\ell_1}(a_{\ell_1})\cap \alpha_{\ell_2}([a_{\ell_2},b_{\ell_2}])=\varnothing$ for $\ell_1\neq\ell_2$, i.e., the endpoints at $\partial B^2$ cannot belong to different curves.

We can now introduce the class of maps we are concerned with.

\begin{definition}\label{def:pL}
	We say that a map $u\in BV(B^2;\gR^2)$ is piecewise Lipschitz if there is a Lipschitz partition $\{\Om_j,j=1,\dots,N\}$ of $B^2$ whose interface is a network such that {the restriction $u\res \Om_j$ belongs to $\text{{\rm Lip}}(\Om_j;\gR^2)$ for every $j$.}
\end{definition}

If $p_i\in B^2$ is a junction point (see (n3)), then it is easily seen that for all $j\in \{1,\dots,N\}$ with $p_i\in \partial \Om_j$ there exists the limit
\begin{equation}\label{limit_corner}
	\beta^i_j:=\lim_{\substack{x\rightarrow p_i\\ x\in \Om_j}}u(x).
\end{equation}

Let $i\in \{1,\dots,m\}$
{be such that $p_i$}
%, suppose that $p_i\in \{p_1,\dots,p_m\}$
is an endpoint of $\alpha_\ell$; by regularity of $\alpha_\ell$ for $\rho>0$ small enough $\alpha_\ell([a_\ell,b_\ell])\cap \partial B_\rho(p_i)$ consists either of a single point, or of two points in the case $\alpha_\ell(a_\ell)= \alpha_\ell(b_\ell)=p_i$.
It follows that the map $u\res  \partial B_\rho(p_i)$ is piecewise Lipschitz and might jump only on  $\Sigma \cap  \partial B_\rho(p_i)$. Hence, the number of these jump
points is, by definition of junction point, $N_i$.
For all such $i$ we denote by $\Om_1^i,\dots,\Om_{N_i}^i,$
the connected components of $B^2\setminus \Sigma$ whose closure contains $p_i$, chosen in counter-clockwise order around $p_i$. By Lipschitz regularity of the partition, any $\Om_j^i$ has a corner at $p_i$ whose aperture is a positive angle $\theta_{j}^i\in(0,2\pi)$. For simplicity of notation denote  $\beta_j^i$ the limit in \eqref{limit_corner} done on $\Om_j^i$ instead of $\Om_j$. For all $i=1,\dots,N_i$ we denote by $\gamma_i$ a Lipschitz curve which parameterizes on $\mathbb S^1$ the polygon  with vertices $\beta_1^i,\beta_2^i,\dots ,\beta_{N_i}^i$ in the order. In \cite{BSS24} the following result is proved:

\begin{theorem}\label{teo1.3}
	Let $u\in BV(B^2;\gR^2)$ be a piecewise Lipschitz map. Then
	\begin{equation}\label{1.2}
		\widetilde {\mathcal A}_{BV}(u,B^2)=A(u,B^2)+\sum_{\ell=1}^n\int_{(a_\ell,b_\ell)\times (0,1)}|\partial_t X_\ell\wedge\partial_sX_\ell| \,\dd t \,\dd s+\sum_{i=1}^m P(\gamma_i),
	\end{equation}
	where
	$X_\ell(t,s)=(t,su^+(\alpha_\ell(t))+(1-s)u^-(\alpha_\ell(t)))$ for $(t,s)\in (a_\ell,b_\ell)\times (0,1)$, with $u^\pm$ being the two traces of $u$ on the jump set $\alpha_\ell((a_\ell,b_\ell))$.
\end{theorem}

In the special case that the function $u$ has no junction points $p_i$, the quantity in \eqref{1.2} (hence with $m=0$ and so that $\sum_{i=1}^m P(\gamma_i)=0$) coincides with
\beq\label{BSS}
\GM(\ol T)=A(u,B^2)+\GM(S_u)
\eeq
for the unique Cartesian current $T_u=G_u+S_u$ in the class $\cT_u^{BV}$.
Therefore, by Theorem \ref{T-cont} we infer that for every continuous integrand $f$ satisfying properties
{$(a),\,(b),\,(c)$ in Sec. \ref{Sec:poly},}
%$(a),\,(b),\,(c)$
we have
\beq\label{BBS2}
\wid\cF_{BV}(u,B^2)=\cF_f(T_u)=\cE_f(u,B^2)+\cF_f(S_u)\,.
\eeq

The previous formula continues to hold if $u\in BV(B^2;\mathbb R^2)$ is piecewise Lipschitz and is such that at every junction point $p_i$ the curve $\gamma_i$ satisfies hypothesis (P0) and (H). In such a case  at any junction point $p_i$ the surface solution of $P(\gamma_i)$ has area equal to the mass of the current ${S_{\gamma_i}}$, {compare Theorem \ref{teo2.9}}.
We have the following:
\begin{theorem}
	Let $u\in BV(B^2;\gR^2)$ be a piecewise Lipschitz map as in Theorem \ref{teo1.3}. Suppose that for all $i=1,\dots,m$ the curve $\gamma_i$ satisfies hypothesis (H), so there is a unique element $(k_1^i,\dots,k^i_{N_i})$ in $\text{{\rm Ad}}_{\min}(\gamma_i)$. Then
	\begin{equation}
		\widetilde {\mathcal F}_{BV}(u,B^2)=\mathcal E_f(u,B^2)+\mathcal F_f(S_u),
	\end{equation}
	where $S_u=T_u-G_u$, and $T_u$ is the unique element in $\mathcal T_u^{BV}$. Explicitly
	\begin{equation}
		\mathcal F_f(S_u)=\sum_{\ell=1}^n\int_{(a_\ell,b_\ell)\times (0,1)}F_f(\alpha_\ell(t),y(s),\tau(t)\wedge w(t))|\partial_t X_\ell\wedge\partial_sX_\ell| \, \dd t\,\dd s +\sum_{i=1}^m P_f(p_i,\gamma_i)
	\end{equation}
	where $y(t)=su^+(t)+(1-s)u^-(t)$, $\tau(t)=\dot\alpha_\ell(t)$, $w(t)=\frac{u^+(t)-u^-(t)}{|u^+(t)-u^-(t)|}$, and for  all $i=1,\dots,m$ it holds
	\begin{equation}
		P_f(p_i,\gamma_i)=\sum_{j=1}^{N_i}k^i_j\int_{U^i_j}F_f(p_i,y,\varepsilon_1\wedge\varepsilon_2)\dd y.
	\end{equation}
\end{theorem}
We recall {that} $X_\ell$ for $(t,s)\in (a_\ell,b_\ell)\times (0,1)$ {is} defined as in Theorem \ref{teo1.3}.

\section{Open questions}\label{sec:open}

We collect some open questions.
\subsection{More general 0-homogeneous maps}
An important question arises if one {wishes} to extend Theorem \ref{teo_0hom1} to the cases in which the {Lipschitz} curve $\gamma$ does not satisfy hypothesis (H).

This is for example the case of the double eight curve $\varphi_8:\mathbb S^1\rightarrow \gR^2$ considered in Proposition \ref{prop:8}. From Remark \ref{rem_28} the corresponding {representative} in $F(2)$ is $\sigma_1\sigma_2\sigma_1^{-1}\sigma_2^{-1}$, and
$\text{{\rm Ad}}_{\min}(\varphi_8)=\{(2,0),(0,2)\}$.
We expect that in such a case the energy gap is 
$$\min \left\{ \int_{C_i}F_f(0_{\gR^2},y,\varepsilon_1\wedge\varepsilon_2)\,\dd y,\;i=1,2 \right\},$$
where $C_1,C_2$ are the unit {disks} with centers $(1,0)$ and $(-1,0)$, respectively. Notice that this is coherent with the case of the area functional. 

\subsection{More on the vortex map: modified catenoid}
Returning {to the vortex map $u_V$ in case of dimension $n=2$,} it would be interesting to obtain for $R>0$ small the optimal estimate in Proposition \ref{Pvortex2}.

Assume for simplicity that $f=f(G)$. In that case, we expect that the energy gap of the {$L^1$} relaxed functional is equal to the $\cF_f$-energy of a suitable vertical current $S_T\in\cR_2(B^2_R\tim\gR^2)$ that is supported on a surface $\Sigma_R$ whose projection onto the domain lies in a segment $I_R$ connecting the origin $0_{\gR^2}$ to a given point $P_R$ at the boundary $\pa B_R$. In strict similarity with the case of the area functional (addressed in \cite{BES1,BES2,BES3}), such a surface $\Sigma_R$ is obtained through a minimization problem for the $\cF_f$-energy in a $3$-dimensional Plateau type problem in $I_R\tim\gR^2$ with a free boundary in correspondence to {the boundary point} $P_R$. In order to apply our continuity results (Proposition \ref{Pcont}) a sequence $v_h\in C^1(B^2_R\times\gR^2)$ must be found such that $G_{v_h}$ converges in mass to the corresponding Cartesian current. In some particular cases in which the function $f$ is explicit, under suitable changes of variables, the aforementioned Plateau type problem with partial free boundary can be rephrased in term of the area functional, but in different domains. In such a case, results in \cite{BMS} will help to study existence and regularity of the corresponding solutions.

\subsection{On the triple junction map}

Considering the relaxed area of the triple junction map $u_T$ introduced in \eqref{triplemap} on {$B^2_R$,} one might ask if something is possible to say about the relaxation in $L^1$ of a more general polyconvex integrand. Following \cite{BePa,Sc20}, the optimal Cartesian current $T\in \cart(B^2_R\times \gR^2)$ obtained as limit of graphs of a recovery sequence for $\widetilde {\mathcal A}_{L^1}(u_T,B^2_R)$ can be written as
$$T=G_{u_T}+S_1+S_2+S_3,$$
where $S_i$ are vertical and live respectively on the three radii in $B^2_R$ forming the jump set $J_{u_T}$ of $u_T$. More precisely, if $R_1$ is the radius of $B^2_R$ corresponding to the polar coordinate $\theta=\pi/3$, $S_1$ is the integration over a surface $\Sigma_1\subset R_1\times \gR^2$ defined as follows: we use the coordinate $\rho\in [0,R)$ to parameterize $R_1$, and choose a Cartesian system of $\gR^2$ such that $\alpha$ and $\beta$ have coordinates $(-\frac12,0)$ and $(\frac12,0)$ respectively. Then we let $Q:=(0,R)\times (-\frac12,\frac12)$ and consider the following Plateau type problem
\begin{equation}
\inf\left\{\int_Q\sqrt{1+|\nabla \psi|^2}\dd x:\psi\in W^{1,1}(Q),\,\psi(\cdot,-\frac12)=\psi(\cdot,\frac12)=0,\,\psi(0,\cdot)=\varphi(\cdot)\right\}
\end{equation}
 where $\varphi:(-\frac12,\frac12)$ is the piecewise affine interpolation map such that $\varphi(-\frac12)=\varphi(\frac12)=0$ and $\varphi(0)=\frac{\sqrt3}{6}$. The resulting minimizer $\psi_{\text{{\rm min}}}$ is a $C^1$ map whose graph is a minimal surface with free boundary on $\{R\}\times \gR^2$; then
 $$\Sigma_1:=G_{\psi_{\text{{\rm min}}}}.$$
The currents $S_2$ and $S_3$ are built in a similar manner.

One might wonder if some analogous argument used to prove Theorem \ref{T-CM} could apply to show that, given a more general polyconvex integrand $f$, the expression of $	\widetilde {\mathcal F}_{L^1}(u_T,B^2_R)$ can be shown. Even more, one may ask if under suitable hypotheses on $f$ there holds
$$	\widetilde {\mathcal F}_{L^1}(u_T,B^2_R)=\mathcal E_f(u_T,B^2_R)+\sum_{i=1}^3\mathcal F_f(S_i).$$

{\bf Acknowledgements.}
D.M. and R.S. are members of the Italian group of Analysis, Probability
and Applications (GNAMPA-INDAM).
\end{document}